\documentclass[review,authoryear]{elsarticle}

\usepackage{amssymb,amsthm,amsmath,graphicx,color}
\usepackage{amsfonts}
\usepackage[caption=false]{subfig}
\usepackage{epsfig,latexsym,graphicx}
\RequirePackage[OT1]{fontenc}
\RequirePackage{amsthm,amsmath}
\RequirePackage{natbib}

\newtheorem{theorem}{Theorem}[section]

\newtheorem{corollary}[theorem]{Corollary}

\theoremstyle{plain}

\begin{document}
	
	\begin{frontmatter}

		\title{Power and Level Robustness of A Composite Hypothesis Testing under Independent Non-Homogeneous Data}
		\author{Abhik Ghosh\fnref{label2}}
		\ead{abhianik@gmail.com}
		\author{Ayanendranath Basu\fnref{label1}}
		\ead{ayanbasu@isical.ac.in}
		\address{Indian Statistical Institute}
		\fntext[label2]{This is a part of the Ph.D. dissertation of the first author.}
		\fntext[label1]{Corresponding Author}

\begin{abstract}
Robust tests of general composite hypothesis under non-identically distributed observations 
is always a challenge. Ghosh and Basu (2018, Statistica Sinica, 28, 1133--1155)
have proposed a new class of test statistics for such problems 
based on the density power divergence, but their robustness with respect to the size and power are not studied in detail.
This note fills this gap by providing a rigorous derivation of 
power and level influence functions of these tests to theoretically justify their robustness.
Applications to the fixed-carrier linear regression model are also provided with empirical illustrations.
\end{abstract}

		\begin{keyword}
{Power Influence Function} \sep 
{Level Influence Function} \sep 
{Robust Hypothesis Testing} \sep 
{Non-Homogeneous Observation} \sep  
{Linear Regression}.
			
		\end{keyword}
		
	\end{frontmatter}

\section{Introduction and Background}\label{SEC:intro}

Robust statistical inference based on non-homogeneous data is always a big challenge
and the likelihood ratio test (LRT), the canonical tool in these situations, is highly sensitive in the presence of outliers.
Literature of alternative robust tests for statistical hypotheses are limited beyond the identically distributed data, 
except for some particular cases like the fixed-carrier linear regression model, etc.
Recently, \cite{Ghosh/Basu:2017} have developed a class of robust testing procedures
under the general set-up of independent but non-homogeneous (INH) observations
based on the  robust estimator of \cite{Ghosh/Basu:2013}.

Under the general INH set-up, we assume that the observations $Y_1, \ldots, Y_n$ are independent but 
$Y_i \sim g_i$ for each $i$ where $g_1, \ldots, g_n$ are potentially different densities with respect 
to some common dominating measure. A parametric family of densities  
${\mathcal F}_{i, {\boldsymbol{\theta}}} = \{f_i(\cdot;{\boldsymbol{\theta}}) |~ {\boldsymbol{\theta}} \in \Theta \}$ 
is assumed to model $g_i$, for each $i=1,2, \ldots, n$, 
and our interest is to make inference about the common parameter  $\boldsymbol{\theta}$.
The most common application is the fixed-carrier regressions, 
where each $f_i(\cdot;{\boldsymbol{\theta}})$ is the (conditional) density of the response 
given the $i$-th (fixed) value of the  covariates. 
In general, we denote by $G_i$ and $F_i(\cdot,{\boldsymbol{\theta}})$ the distribution functions of $g_i$ 
and $f_i(\cdot;{\boldsymbol{\theta}})$ respectively. 
Under this INH set-up, \cite{Ghosh/Basu:2013} have developed a general robust estimator of $\boldsymbol{\theta}$
using the density power divergence (DPD) of \cite{Basu/etc:1998}; 
this DPD measure, having a tuning parameter $\tau$, is defined between densities $f_1$ and $f_2$ as
\begin{equation}\label{EQ:dpd}
d_\tau(f_1,f_2) = \displaystyle \left\{\begin{array}{ll}
\displaystyle \int  \left[f_2^{1+\tau} - \left(1 + \frac{1}{\tau}\right)  f_2^\tau f_1 + 
\frac{1}{\tau} f_1^{1+\tau}\right], & {\rm for} ~\tau > 0,\\
\displaystyle \int f_1 \log(f_1/f_2), & {\rm for} ~\tau = 0.  
\end{array}\right.
\end{equation}
Since there are $n$ different densities for INH set-up,
\cite{Ghosh/Basu:2013} minimized the average DPD measure 
$\frac{1}{n} \sum_{i=1}^n d_\tau(\widehat{g}_i(.),f_i(.;{\boldsymbol{\theta}}))$ 
with respect to $\boldsymbol{\theta}\in \Theta$, 
where $\widehat{g}_i$ is an estimator of $g_i$ based on the empirical distribution function. 
This minimum DPD estimator (MDPDE) has high efficiency and robustness properties, controlled by $\tau$,
and works well in different fixed-design regressions by 
\cite{Ghosh/Basu:2013, Ghosh/Basu:2014} and \cite{Ghosh:2017a,Ghosh:2017b}.
At $\tau=0$, the MDPDE coincides with the maximum likelihood estimator (MLE).
Using this MDPDE, \cite{Ghosh/Basu:2017} have developed a class of robust DPD based tests 
for both simple and composite hypotheses indexed by the same $\tau$;
they coincide with the LRT at $\tau=0$ and 
provide its robust generalization at $\tau>0$ without significant loss in efficiency. 
However, their theoretical robustness properties need to be studied in greater detail,
particularly for composite hypothesis testing problems,
where no details about the size and power robustness are available.

Since the size and power are the two most important measures to study the performance of any test,
in this paper, we present detailed analysis for such robustness issues for the composite hypothesis tests
of \cite{Ghosh/Basu:2017}. In particular, we study their power and level influence functions 
to justify their robustness with a concrete theory;
this needs some non-trivial extensions of the corresponding results from simple hypothesis case.
We also illustrate their applications in testing general linear hypothesis 
under a fixed-carrier linear regression model (LRM) with unknown error variance. 
Empirical results from an extensive simulation study second our theoretical robustness analyses.

We provide a brief description of 
the composite hypothesis tests from \cite{Ghosh/Basu:2017} in Section \ref{SEC:8composite_testing}.
Our main results about the level and power influence functions are provided in Section \ref{SEC:8IF_power_composite}.
Section \ref{SEC:linear_regression} presents the application to the LRMs
and numerical illustrations are given in Section \ref{SEC:9simulation_test}.
Concluding remarks are given in Section \ref{SEC:conclusion}. 
All notations are given in \ref{APP:A}, 
whereas the required assumptions and some background results 
are presented in the Online Supplement for completeness.

\section{DPD based  Tests for Composite Hypotheses under the INH Set-up}\label{SEC:8composite_testing}

Consider the INH set-up of Section \ref{SEC:intro} and 
the problem of testing the composite hypothesis of the form
\begin{equation}\label{EQ:8composite_hypo}
 H_0 : {\boldsymbol{\theta}} \in \Theta_0 ~~~\mbox{ against }~~~~ H_1 : {\boldsymbol{\theta}}  \notin \Theta_0,
\end{equation}
where $\Theta_0\subset\Theta$. In most applications, 
the (fixed) null parameter space $\Theta_0$ is defined in terms of $r$ independent restrictions, say
$\boldsymbol\upsilon({\boldsymbol{\theta}}) = \boldsymbol{0}_r$. 
\cite{Ghosh/Basu:2017} have proposed to test (\ref{EQ:8composite_hypo}) by 
the DPD based test statistics 

\begin{eqnarray}\label{EQ:8DPDTS2}
S_{\gamma}( {{\boldsymbol{\theta}}_n^\tau}, {\widetilde{{\boldsymbol{\theta}}}_n^\tau}) 
= 2 \sum_{i=1}^n d_\gamma(f_i(.;{\boldsymbol{\theta}}_n^\tau),f_i(.;\widetilde{{\boldsymbol{\theta}}}_n^\tau)),
~~~\tau,\gamma \geq 0,
\end{eqnarray}
where ${\boldsymbol{\theta}}_n^\tau$ and $\widetilde{{\boldsymbol{\theta}}}_n^\tau$ 
are the MDPDE and the restricted MDPDE (RMDPDE) of $\boldsymbol{\theta}$ respectively;
the RMDPDE has to be obtained by minimizing the average DPD measure only over $\boldsymbol{\theta}\in \Theta_0$
(See Results 1 and 2 in Online Supplement for their asymptotic distributions).
\cite{Ghosh/Basu:2017} have shown that, in general, its asymptotic null distribution 
is a linear combination of (central) chi-square distributions (Result 3 in Online Supplement);
some suitable approximations are also suggested for its critical values following  \cite{Basu/etc:2013a}. 
Further, this DPD based test is consistent at any fixed alternative. 

However, in terms of robustness, only the influence function (IF) of the test statistic 
have been discussed in \cite{Ghosh/Basu:2017}.
The statistical functional corresponding to the test statistics in (\ref{EQ:8DPDTS2}) is defined as
\begin{eqnarray}\label{EQ:8DPDTS2_IF}
S_{\gamma, \tau}(\underline{\mathbf{G}}) = \sum_{i=1}^n 
d_\gamma(f_i(.;U_\tau(\underline{\mathbf{G}})),f_i(.;\widetilde{U}_\tau(\underline{\mathbf{G}}))),
\nonumber
\end{eqnarray}
where $\underline{\mathbf{G}} = (G_1, \cdots, G_n)$
and $U_\tau(\underline{\mathbf{G}})$ and $\widetilde{U}_\tau(\underline{\mathbf{G}})$ are the functionals 
corresponding to the MDPDE and the RMDPDE, respectively,  defined as the minimizers of 
$\frac{1}{n} \sum_{i=1}^n d_\tau({g}_i(.),f_i(.;{\boldsymbol{\theta}}))$ 
with respect to $\boldsymbol{\theta}\in \Theta$ and $\boldsymbol{\theta}\in \Theta_0$.
Consider contamination in all densities at the contamination points 
in $\mathbf{t} = (t_1, \ldots, t_n)$ respectively.
When evaluating at the null distribution $\underline{\mathbf{G}}=\underline{\mathbf{F}}_{\theta_0}
= \left(F_1(\cdot, \boldsymbol{\theta}_0), \ldots, F_n(\cdot, \boldsymbol{\theta}_0)\right)$ 
with $\theta_0 \in \Theta_0$, the first order IF of $S_{\gamma, \tau}$ is  identically zero and 
the corresponding second order IF is \citep{Ghosh/Basu:2017}
\begin{eqnarray}
\mathcal{IF}_2(\mathbf{t}, S_{\gamma, \tau}, \underline{\mathbf{F}}_{\theta_0}) 
= n \cdot \boldsymbol{D}_{\tau}(\mathbf{t}, \boldsymbol{\theta}_0)^T \boldsymbol{A}_n^\gamma(\boldsymbol{\theta}_0) 
\boldsymbol{D}_{\tau}(\mathbf{t}, \boldsymbol{\theta}_0),
\label{EQ:IF2_TS}
\end{eqnarray}
where $\boldsymbol{D}_{\tau}(\mathbf{t}, \boldsymbol{\theta}_0)
= \left[\mathcal{IF}(\mathbf{t}, \boldsymbol{U}_\tau, \underline{\mathbf{F}}_{\theta_0}) - 
\mathcal{IF}(\mathbf{t}, \widetilde{\boldsymbol{U}}_\tau, \underline{\mathbf{F}}_{\theta_0})\right]$,
the difference between IFs  of the MDPDE $\boldsymbol{U}_\tau$ and 
the RMDPDE $\widetilde{\boldsymbol{U}}_\tau$ at $\underline{\mathbf{F}}_{\theta_0}$.
But, both these IFs are both bounded at $\tau>0$ for most parametric models; 
at $\tau=0$ the IF of the MDPDE (MLE) is unbounded but that of RMDPDE depends 
on the restrictions $\boldsymbol{\upsilon}=\boldsymbol{0}$.
So the second order IF (\ref{EQ:IF2_TS}) of our test statistics is bounded 
whenever 	$\boldsymbol{D}_{\tau}(\mathbf{t}, \boldsymbol{\theta}_0)$ is bounded, 
i.e., the IFs of MDPDE and RMDPDE both are bounded or both diverge at the same rate;
this holds for $\tau>0$ in most cases. 
At $\tau=0$, this new test coincides with the non-robust LRT having unbounded IF.

\section{Power and Level Influence Functions}
\label{SEC:8IF_power_composite}

For a hypothesis testing procedure, it is not enough to study only the properties of the test statistics;
the level and power are two basic components of hypothesis testing whose robustness is essential 
to fully justify a new robust test procedure. In this section,
we study the theoretical robustness properties of the power and level of 
the DPD based test in (\ref{EQ:8DPDTS2}); 
it is done through the examination of classical power influence functions (PIF)
and level influence function (LIF). 

The PIF and LIF of a test measure the effect of infinitesimal contamination on its power and level respectively.  
However, the DPD based test (\ref{EQ:8DPDTS2}) is consistent at any fixed alternative \citep{Ghosh/Basu:2017}
and hence its power against any fixed alternative is always one.
Further, exact finite-sample power is much difficult to derive. 
So, we study the effect of contamination on its asymptotic power against a sequence of contiguous alternatives 
$H_{1,n} : {\boldsymbol{\theta}}={\boldsymbol{\theta}}_n$, where ${\boldsymbol{\theta}}_n = {\boldsymbol{\theta}}_0 + n^{-1/2}\boldsymbol\Delta$ 
with ${\boldsymbol{\theta}}_0\in \Theta_0$ and $\boldsymbol\Delta \in \mathbb{R}^p - \{\boldsymbol{0}_p\}$. 
Such a ${\boldsymbol{\theta}}_0$ must be a limit point of $\Theta_0$;
we assume $\Theta_0$ to be closed ensuring the existence of such a sequence $\boldsymbol{\theta}_n \in \Theta$. 
Then, we consider the contamination over these contiguous alternatives in such a way that
the contamination effect vanishes at the same rate as ${\boldsymbol{\theta}}_n\rightarrow{\boldsymbol{\theta}}_0$ 
when $n\rightarrow\infty$; 
this is necessary to make the neighborhood of the null and alternative hypotheses well separated \citep{Hampel/etc:1986}. 
Note that $\boldsymbol\Delta = \boldsymbol{0}$ yields the results associated with level of the test.
Thus, assuming contamination in all densities as in the previous section, 
the contaminated distributions need to be defined as
$$
\mathbf{\underline{F}}_{n,\epsilon,\mathbf{t}}^P = \left(1-\frac{\epsilon}{\sqrt{n}}\right) 
\mathbf{\underline{F}}_{{\boldsymbol{\theta}}_n}+\frac{\epsilon}{\sqrt{n}} \wedge_{\mathbf{t}},
\mbox{ and }
\mathbf{\underline{F}}_{n,\epsilon,\mathbf{t}}^L = \left(1-\frac{\epsilon}{\sqrt{n}}\right) 
\mathbf{\underline{F}}_{{\boldsymbol{\theta}}_0}+\frac{\epsilon}{\sqrt{n}} \wedge_{\mathbf{t}},
$$
for studying the stability of power and level respectively,
where $\epsilon$ is the contamination proportion
and $ \wedge_{\mathbf{t}}=( \wedge_{{t}_1}, \ldots,  \wedge_{{t}_n})$ with $ \wedge_{{t}_i}$
being the degenerate distribution at $t_i$ for each $i=1, \ldots, n$.
Then the PIF and LIF of the test in (\ref{EQ:8DPDTS2}), at the significance level $\alpha$, 
are defined, see \cite{Hampel/etc:1986}, as
\begin{eqnarray}
PIF(\mathbf{t}; S_{\gamma, \tau}, \mathbf{\underline{F}}_{{\boldsymbol{\theta}}_0} ) 
&=& \lim_{n \rightarrow \infty} ~ \frac{\partial}{\partial \epsilon} 
P_{\mathbf{\underline{F}}_{n,\epsilon,\mathbf{t}}^P }( S_{\gamma}( {{\boldsymbol{\theta}}_n^\tau}, \widetilde{{\boldsymbol{\theta}}}_n^\tau) 
>  s_\alpha^{\tau,\gamma}) \big|_{\epsilon=0},\nonumber\\
LIF(\mathbf{t}; S_{\gamma, \tau}, \mathbf{\underline{F}}_{{\boldsymbol{\theta}}_0} ) 
&=& \lim_{n \rightarrow \infty} ~ \frac{\partial}{\partial \epsilon} 
P_{\mathbf{\underline{F}}_{n,\epsilon,\mathbf{t}}^L }
(S_{\gamma}( {{\boldsymbol{\theta}}_n^\tau}, \widetilde{{\boldsymbol{\theta}}}_n^\tau)> s_\alpha^{\tau,\gamma}) 
\big|_{\epsilon=0},
\nonumber
\end{eqnarray}
where $s_\alpha^{\tau,\gamma}$ is the $(1-\alpha)$-th quantile of the asymptotic null distribution of
$S_{\gamma}( {{\boldsymbol{\theta}}_n^\tau}, \widetilde{{\boldsymbol{\theta}}}_n^\tau)$.
\cite{Ghosh/Basu:2017} have discussed these LIF and PIF for testing the simple null hypothesis;
further applications can be found in 
\cite{Huber/Carol:1970}, \cite{Heritier/Ronchetti:1994} and \cite{Toma/Broniatowski:2010}
for both types of hypotheses. Following the same line of arguments, 
we start with the derivation of the asymptotic power of the DPD based test (\ref{EQ:8DPDTS2})
under $\mathbf{\underline{F}}_{n,\epsilon,\mathbf{y}}^P$, recalling the notations from \ref{APP:A}.

\begin{theorem}
Suppose that Assumptions (A1)-(A10), given in Online Supplement, 
hold at $\boldsymbol\theta=\boldsymbol\theta_0$ under the INH set-up. 
Then, for any $\boldsymbol\Delta \in \mathbb{R}^p$, $\epsilon \geq 0$, we have the following results.
	\begin{enumerate}
		\item[(i)] Under $\mathbf{\underline{F}}_{n,\epsilon,\mathbf{t}}^P $, 
		$S_{\gamma}( {{\boldsymbol{\theta}}_n^\tau}, {\widetilde{{\boldsymbol{\theta}}}_n^\tau}) 
		\displaystyle\mathop{\rightarrow}^\mathcal{D} \boldsymbol{W}^T  \boldsymbol{A}_\gamma(\theta_0)\boldsymbol{W}$,
		where $\boldsymbol{W}\sim N_p\left(\widetilde{\boldsymbol{\Delta}^*}, 
		\widetilde{\boldsymbol{\Sigma}}_\tau(\boldsymbol{\theta}_0)\right)$ with 		
		$\widetilde{\boldsymbol{\Delta}^*} = \left[\boldsymbol{\Delta} 
		+ \epsilon \boldsymbol{D}_{\tau}(\mathbf{t}, \boldsymbol{\theta}_0)\right]$.
		
		\item[(ii)] Suppose the $r$ eigenvalues of 
		${\boldsymbol{A}_{\gamma}}({\boldsymbol{\theta}}_0)\widetilde{\boldsymbol\Sigma}_\tau({\boldsymbol{\theta}}_0)$ 
		are denoted as $\widetilde{\zeta_1^{\gamma, \tau}}({\boldsymbol{\theta}}_0)$, $\ldots$, 
		$\widetilde{\zeta_r^{\gamma, \tau}}({\boldsymbol{\theta}}_0)$ 
		with the corresponding  normalized eigenvector matrix being
		$\widetilde{\boldsymbol{P}}_{\tau,\gamma}(\theta_0)$. Denote 
		$
		\widetilde{\boldsymbol{P}}_{\tau,\gamma}(\boldsymbol{\theta}_0)
		\widetilde{\boldsymbol{\Sigma}}_\tau^{-1/2}(\boldsymbol{\theta}_0)
		\widetilde{\boldsymbol\Delta^*} =\left({\widetilde{\delta}_1}, \ldots, {\widetilde{\delta}_p}\right)^T. 
		$
		Then, the asymptotic distribution in (i) is also the distribution of 
		$\sum\limits_{i=1}^{r}\widetilde{\zeta_i^{\gamma, \tau}}(\theta_0)\chi_{1,\widetilde{\delta}_i}^2$,
		where $\chi_{1,\widetilde{\delta}_i^2}^2$s are independent non-central chi-square variables with 
		degrees of freedom (df) $1$ and non-centrality parameter (ncp) $\widetilde{\delta}_i^2$ respectively, 
		for $i=1, \ldots, r$.
		
		\item[(iii)] The asymptotic power of the DPD based test (\ref{EQ:8DPDTS2}) under 
		$\mathbf{\underline{F}}_{n,\epsilon,\mathbf{t}}^P$ is given by
		\begin{eqnarray}
			P_{\tau,\gamma}^*(\boldsymbol\Delta, \epsilon; \alpha) &=& 
			\lim_{n \rightarrow \infty}  P_{\mathbf{\underline{F}}_{n,\epsilon,\mathbf{t}}^P }\left( 
			S_{\gamma}( {\theta_n^\tau}, \widetilde{\theta}_n^\tau) >  s_\alpha^{\tau,\gamma}\right)
=\sum\limits_{v=0}^{\infty}~\widetilde{C_v^{\gamma,\tau}}(\boldsymbol\theta_0, \widetilde{\boldsymbol\Delta^*})
P\left(\chi_{r+2v}^2 > {s_\alpha^{\tau, \gamma}}/{\widetilde{\zeta_{(1)}^{\gamma,\tau}}(\boldsymbol\theta_0)}\right), 
\nonumber
		\end{eqnarray}
	where  	$\widetilde{\zeta_{(1)}^{\gamma,\tau}}(\boldsymbol\theta_0) =
	\min\limits_{1\leq r\leq r}\widetilde{\zeta_{i}^{\gamma,\tau}}(\boldsymbol\theta_0)$,
	$\chi_{r+2v}^2$ are independent chi-squares with df $r+2v$ 	for $v\geq 0$,  and
		$$
		\widetilde{C_v^{\gamma,\tau}}(\boldsymbol\theta_0, \widetilde{\boldsymbol\Delta^*}) = \frac{1}{v!} \left(
		\prod\limits_{j=1}^{r}\frac{\widetilde{\zeta_{(1)}^{\gamma,\tau}}(\boldsymbol\theta_0)}{
			\widetilde{\zeta_{j}^{\gamma,\tau}}(\boldsymbol\theta_0)}\right)^{1/2}
		e^{-\frac{1}{2}\sum\limits_{j=1}^{r}\widetilde{\delta}_j^2} E(\widetilde{Q}^v), \nonumber
		$$
		$$
		\mbox{with }~~~~~~~~~
		\widetilde{Q} = \frac{1}{2}\sum\limits_{j=1}^{r}\left[\left(1 - 
		\frac{\widetilde{\zeta_{(1)}^{\gamma,\tau}}(\boldsymbol\theta_0)}{
			\widetilde{\zeta_{j}^{\gamma,\tau}}(\boldsymbol\theta_0)}\right)^{1/2}Z_j
		+ {\widetilde{\delta}_j}\left(\frac{\widetilde{\zeta_{(1)}^{\gamma,\tau}}(\boldsymbol\theta_0)}{
			\widetilde{\zeta_{j}^{\gamma,\tau}}(\boldsymbol\theta_0)}\right)^{1/2}\right]^2,
		$$
		for $r$ independent standard normal random variables $Z_1, \ldots, Z_r$.
	\end{enumerate}
	\label{THM:8asymp_Contaminated_power_composite}
\end{theorem}
\noindent
\textbf{Proof}
All notations and matrices used in this proof are defined in \ref{APP:A} for brevity.
Let us denote ${\boldsymbol{\theta}}_n^* = \boldsymbol{U}_{\tau}(\mathbf{\underline{F}}_{n,\epsilon,\mathbf{t}}^P)$ 
and $\widetilde{{\boldsymbol{\theta}}_n^*} = \widetilde{\boldsymbol{U}}_{\tau}(\mathbf{\underline{F}}_{n,\epsilon,\mathbf{t}}^P)$. 
Fix any $i=1, \ldots, n$. We consider the second order Taylor series expansion of 
$d_\gamma(f_i(\cdot;{\boldsymbol{\theta}}),f_i(\cdot;\widetilde{{\boldsymbol{\theta}}}_n^\tau))$ 
around ${\boldsymbol{\theta}} = {\boldsymbol{\theta}}_n^*$ at ${\boldsymbol{\theta}} = {\boldsymbol{\theta}}_n^\tau $ as,
\begin{eqnarray}
&&d_\gamma(f_i(\cdot;{\boldsymbol{\theta}}_n^\tau),f_i(\cdot;\widetilde{{\boldsymbol{\theta}}}_n^\tau)) = 
d_\gamma(f_i(\cdot;{\boldsymbol{\theta}}_n^*),f_i(\cdot;\widetilde{{\boldsymbol{\theta}}}_n^\tau)) + 
\boldsymbol{M}_{1, \gamma}^{(i)}({\boldsymbol{\theta}}_n^*, \widetilde{{\boldsymbol{\theta}}}_n^\tau)^T ({\boldsymbol{\theta}}_n^\tau - {\boldsymbol{\theta}}_n^*) 
\nonumber \\ &&~~~~~~~~~~~~~~~~
+\frac{1}{2} ({\boldsymbol{\theta}}_n^\tau - {\boldsymbol{\theta}}_n^*)^T 
\boldsymbol{A}_{1,1,\gamma}^{(i)}({\boldsymbol{\theta}}_n^*, \widetilde{{\boldsymbol{\theta}}}_n^\tau) 
({\boldsymbol{\theta}}_n^\tau - {\boldsymbol{\theta}}_n^*) + o(||{\boldsymbol{\theta}}_n^\tau - {\boldsymbol{\theta}}_n^*||^2).
~~~~~~ 
\label{EQ:8taylor_power_composite}
\end{eqnarray}
Now, using Result 1 of Online Supplement and 
the consistency of ${\boldsymbol{\theta}}_n^*$ we know that,  
under  $\mathbf{\underline{F}}_{n,\epsilon,\mathbf{t}}^P$,
$
\sqrt{n} ({\boldsymbol{\theta}}_n^\tau - {\boldsymbol{\theta}}_n^*) \mathop{\rightarrow}^\mathcal{D}  
N_p(\boldsymbol{0}, \boldsymbol{J}_\tau^{-1}({\boldsymbol{\theta}}_0) \boldsymbol{V}_\tau({\boldsymbol{\theta}}_0) 
\boldsymbol{J}_\tau^{-1}({\boldsymbol{\theta}}_0) ). 
$
Further Taylor series expansions around ${\boldsymbol{\theta}}=\widetilde{{\boldsymbol{\theta}}_n^*}$ at 
${\boldsymbol{\theta}}=\widetilde{{\boldsymbol{\theta}}}_n^\tau$ lead to 
\begin{eqnarray}
&&d_\gamma(f_i(\cdot;{\boldsymbol{\theta}}_n^*),f_i(\cdot;\widetilde{{\boldsymbol{\theta}}}_n^\tau)) = 
d_\gamma(f_i(\cdot;{\boldsymbol{\theta}}_n^*),f_i(\cdot;\widetilde{{\boldsymbol{\theta}}_n^*})) 
+ \boldsymbol{M}_{2, \gamma}^{(i)}({\boldsymbol{\theta}}_n^*, \widetilde{{\boldsymbol{\theta}}_n^*})^T 
(\widetilde{{\boldsymbol{\theta}}}_n^\tau - \widetilde{{\boldsymbol{\theta}}_n^*}) \nonumber \\ 
& &~~~~~~~~~~~~+\frac{1}{2} (\widetilde{{\boldsymbol{\theta}}}_n^\tau - \widetilde{{\boldsymbol{\theta}}_n^*})^T 
\boldsymbol{A}_{2,2,\gamma}^{(i)}({{\boldsymbol{\theta}}}_n^*, \widetilde{{\boldsymbol{\theta}}_n^*}) 
(\widetilde{{\boldsymbol{\theta}}}_n^\tau - \widetilde{{\boldsymbol{\theta}}_n^*}) 
+ o(||\widetilde{{\boldsymbol{\theta}}}_n^\tau - \widetilde{{\boldsymbol{\theta}}_n^*}||^2), \nonumber\\
&&\boldsymbol{M}_{1, \gamma}^{(i)}({\boldsymbol{\theta}}_n^*, \widetilde{{\boldsymbol{\theta}}}_n^\tau) 
= \boldsymbol{M}_{1, \gamma}^{(i)}({\boldsymbol{\theta}}_n^*, \widetilde{{\boldsymbol{\theta}}_n^*}) 
+ \boldsymbol{A}_{2,1,\gamma}^{(i)}({{\boldsymbol{\theta}}}_n^*, \widetilde{{\boldsymbol{\theta}}_n^*}) (\widetilde{{\boldsymbol{\theta}}}_n^\tau - \widetilde{{\boldsymbol{\theta}}_n^*}) 
+ o(||\widetilde{{\boldsymbol{\theta}}}_n^\tau - \widetilde{{\boldsymbol{\theta}}_n^*}||),\nonumber
\end{eqnarray}
and $\boldsymbol{A}_{1,1, \gamma}^{(i)}({\boldsymbol{\theta}}_n^*, \widetilde{{\boldsymbol{\theta}}}_n^\tau) =  
\boldsymbol{A}_{1,1,\gamma}^{(i)}({{\boldsymbol{\theta}}}_n^*, \widetilde{{\boldsymbol{\theta}}_n^*}) + o_P(1)$.
Again, for each $j=1, 2$, Taylor series expansion of 
$\boldsymbol{M}_{j, \gamma}^{(i)}({\boldsymbol{\theta}}, \widetilde{{\boldsymbol{\theta}}_n^*})$ 
around ${\boldsymbol{\theta}} = {\boldsymbol{\theta}}_0$ at ${\boldsymbol{\theta}} = {\boldsymbol{\theta}}_n^*$ gives 
\begin{eqnarray}
\boldsymbol{M}_{j, \gamma}^{(i)}({\boldsymbol{\theta}}_n^*, \widetilde{{\boldsymbol{\theta}}_n^*}) 
&=& \boldsymbol{M}_{j, \gamma}^{(i)}({\boldsymbol{\theta}}_0, \widetilde{{\boldsymbol{\theta}}_n^*}) 
+ \frac{1}{\sqrt{n}}  \boldsymbol{A}_{1,j, \gamma}^{(i)}({\boldsymbol{\theta}}_0, \widetilde{{\boldsymbol{\theta}}_n^*}){\boldsymbol{\Delta}} 
+ \frac{\epsilon}{\sqrt{n}} \boldsymbol{A}_{1,j, \gamma}^{(i)}({\boldsymbol{\theta}}_0, \widetilde{{\boldsymbol{\theta}}_n^*})
\mathcal{IF}(\mathbf{t}; \boldsymbol{U}_{\tau}, \mathbf{\underline{F}}_{{\boldsymbol{\theta}}_0}) 
+ o\left(\frac{1}{\sqrt{n}}\right)
\nonumber\\
&=& \boldsymbol{M}_{j, \gamma}^{(i)}({\boldsymbol{\theta}}_0, \widetilde{{\boldsymbol{\theta}}_n^*}) 
+ \frac{1}{\sqrt{n}}  \boldsymbol{A}_{1,j, \gamma}^{(i)}({\boldsymbol{\theta}}_0, \widetilde{{\boldsymbol{\theta}}_n^*})\widetilde{{\boldsymbol{\Delta}}_1} 
+ o\left(\frac{1}{\sqrt{n}}\right),\nonumber
\end{eqnarray}
where $\widetilde{{\boldsymbol{\Delta}}_1} = {\boldsymbol{\Delta}} 
+ \epsilon \mathcal{IF}(\mathbf{t}; \boldsymbol{U}_{\tau}, \mathbf{\underline{F}}_{{\boldsymbol{\theta}}_0})$.
For each $j,k=1,2$ and $i=1,\ldots,n$, similar use of suitable Taylor series expansions yield
$\boldsymbol{A}_{j,k,\gamma}^{(i)}({{\boldsymbol{\theta}}}_n^*, \widetilde{{\boldsymbol{\theta}}_n^*})
= \boldsymbol{A}_{1,1,\gamma}^{(i)}({{\boldsymbol{\theta}}}_0, \widetilde{{\boldsymbol{\theta}}_n^*}) + o(1)
=\boldsymbol{A}_{1,1,\gamma}^{(i)}({{\boldsymbol{\theta}}}_n^*, {\boldsymbol{\theta}}_0) + o(1)$, and
\begin{eqnarray}
\boldsymbol{M}_{j, \gamma}^{(i)}({\boldsymbol{\theta}}_0, \widetilde{{\boldsymbol{\theta}}_n^*}) 
&=& \boldsymbol{M}_{j, \gamma}^{(i)}({\boldsymbol{\theta}}_0, {{\boldsymbol{\theta}}_0}) + \frac{\epsilon}{\sqrt{n}} 
\boldsymbol{A}_{2,j, \gamma}^{(i)}({\boldsymbol{\theta}}_0, {\boldsymbol{\theta}}_0) 
\mathcal{IF}(\mathbf{t}; \widetilde{\boldsymbol{U}}_\tau, \mathbf{\underline{F}}_{{\boldsymbol{\theta}}_0})
+ o\left(\frac{1}{\sqrt{n}}\right),
\nonumber\\
\boldsymbol{M}_{j, \gamma}^{(i)}({{\boldsymbol{\theta}}_n^*}, {\boldsymbol{\theta}}_0) 
&=& \boldsymbol{M}_{j, \gamma}^{(i)}({\boldsymbol{\theta}}_0, {{\boldsymbol{\theta}}_0}) + \frac{1}{\sqrt{n}} 
\boldsymbol{A}_{1,j, \gamma}^{(i)}({\boldsymbol{\theta}}_0, {\boldsymbol{\theta}}_0) \widetilde{{\boldsymbol{\Delta}}_1} + o\left(\frac{1}{\sqrt{n}}\right).
\nonumber
\end{eqnarray}
Now, we use these expressions to simplify Equation (\ref{EQ:8taylor_power_composite})
and consider its summation over all $i=1, \ldots, n$.
But, we also know that ${{\boldsymbol{\theta}}}_n^*\rightarrow {\boldsymbol{\theta}}_0$ as $n \rightarrow\infty$
and so $\frac{1}{n}\sum\limits_{i=1}^{n} \boldsymbol{M}_{j, \gamma}^{(i)}({{\boldsymbol{\theta}}_n^*}, {\boldsymbol{\theta}}_0) 
\rightarrow \boldsymbol{M}_\gamma({\boldsymbol{\theta}}_0)=\boldsymbol{0}$ and
$\frac{1}{n}\sum\limits_{i=1}^{n} \boldsymbol{A}_{j, k, \gamma}^{(i)}({{\boldsymbol{\theta}}_n^*}, {\boldsymbol{\theta}}_0) 
\rightarrow (-1)^{j+k} \boldsymbol{A}_\gamma({\boldsymbol{\theta}}_0)$ as $n\rightarrow\infty$ for all $j,k =1, 2$. 
Thus, we get 
\begin{eqnarray}
\frac{1}{\sqrt{n}}\sum\limits_{i=1}^{n} \boldsymbol{M}_{j, \gamma}^{(i)}({\boldsymbol{\theta}}_n^*, \widetilde{{\boldsymbol{\theta}}_n^*}) 
&=&  (-1)^{1+j} \boldsymbol{A}_{\gamma}({\boldsymbol{\theta}}_0)\left[\widetilde{{\boldsymbol{\Delta}}_1} - 
\epsilon \mathcal{IF}(\mathbf{t}; \widetilde{\boldsymbol{U}}_\tau, \mathbf{\underline{F}}_{{\boldsymbol{\theta}}_0})\right]+ o\left(1\right)
=(-1)^{1+j} \boldsymbol{A}_{\gamma}({\boldsymbol{\theta}}_0)\widetilde{{\boldsymbol{\Delta}}^*} + o\left(1\right),
\nonumber
\end{eqnarray}
where $\widetilde{{\boldsymbol{\Delta}}^*}$ is as defined in the theorem. 
Next, another Taylor series expansion of $d_\gamma(f_i(\cdot;{\boldsymbol{\theta}}_n^*),f_i(\cdot;{\boldsymbol{\theta}}))$ 
around ${\boldsymbol{\theta}} = {\boldsymbol{\theta}}_0$ at ${\boldsymbol{\theta}} = \widetilde{{\boldsymbol{\theta}}_n^*} $ gives
\begin{eqnarray}
d_\gamma(f_i(\cdot;{\boldsymbol{\theta}}_n^*), f_i(\cdot;\widetilde{{\boldsymbol{\theta}}_n^*})) &=& 
d_\gamma(f_i(\cdot;{\boldsymbol{\theta}}_n^*),f_i(\cdot;{{\boldsymbol{\theta}}}_0)) + \frac{\epsilon}{\sqrt{n}}
\boldsymbol{M}_{2, \gamma}^{(i)}({\boldsymbol{\theta}}_n^*, {\boldsymbol{\theta}}_0)
\mathcal{IF}(\mathbf{t}; \widetilde{\boldsymbol{U}}_\tau, \mathbf{\underline{F}}_{{\boldsymbol{\theta}}_0})  
\nonumber \\ 
&&+\frac{\epsilon^2}{2n} 
\mathcal{IF}(\mathbf{t}; \widetilde{\boldsymbol{U}}_\tau, \mathbf{\underline{F}}_{{\boldsymbol{\theta}}_0})^T 
\boldsymbol{A}_{2,2,\gamma}^{(i)}({\boldsymbol{\theta}}_n^*, {\boldsymbol{\theta}}_0)
\mathcal{IF}(\mathbf{t}; \widetilde{\boldsymbol{U}}_\tau, \mathbf{\underline{F}}_{{\boldsymbol{\theta}}_0}) 
+ o\left(\frac{1}{n}\right). \nonumber \\
\mbox{Similarly }~~
2 \sum_{i=1}^n d_{\gamma}(f_i(\cdot;{\boldsymbol{\theta}}_n^*), f_i(\cdot;{\boldsymbol{\theta}}_0)) &=& 
{\boldsymbol{\Delta}}^T  \boldsymbol{A}_n^\gamma({\boldsymbol{\theta}}_0){\boldsymbol{\Delta}} + 2\epsilon {\boldsymbol{\Delta}}^T \boldsymbol{A}_n^\gamma({\boldsymbol{\theta}}_0) 
\mathcal{IF}(\mathbf{t}; \boldsymbol{U}_\tau, \mathbf{\underline{F}}_{{\boldsymbol{\theta}}_0}) 
\nonumber \\&+& 
\epsilon^2   
\mathcal{IF}(\mathbf{t}; \boldsymbol{U}_\tau, \mathbf{\underline{F}}_{{\boldsymbol{\theta}}_0}) 
\boldsymbol{A}_n^\gamma({\boldsymbol{\theta}}_0)^T 
\mathcal{IF}(\mathbf{t}; \boldsymbol{U}_\tau, \mathbf{\underline{F}}_{{\boldsymbol{\theta}}_0}) + o\left(1\right) \nonumber \\
&=& \widetilde{{\boldsymbol{\Delta}}_1}^T  \boldsymbol{A}_n^\gamma({\boldsymbol{\theta}}_0)\widetilde{{\boldsymbol{\Delta}}_1} + o\left(1\right).\nonumber
\end{eqnarray}
Combining last two equations,
$2 \sum_{i=1}^n d_{\gamma}(f_i(\cdot;{\boldsymbol{\theta}}_n^*), f_i(\cdot;\widetilde{{\boldsymbol{\theta}}_n^*})) 
=\widetilde{{\boldsymbol{\Delta}}^*}^T  \boldsymbol{A}_\gamma({\boldsymbol{\theta}}_0)
\widetilde{{\boldsymbol{\Delta}}^*} + o\left(1\right)$.
%
Therefore, noting that $n \times o(||{\boldsymbol{\theta}}_n^\tau - {\boldsymbol{\theta}}_n^*||^2) = o_P(1)$ and 
$n \times o(||\widetilde{{\boldsymbol{\theta}}}_n^\tau - \widetilde{{\boldsymbol{\theta}}_n^*}||^2) = o_P(1)$, 
we get the simplified expression as follows. 
\begin{eqnarray}
&&2\sum_{i=1}^n d_{\gamma}(f_i(\cdot;{\boldsymbol{\theta}}_n^\tau), f_i(\cdot;\widetilde{{\boldsymbol{\theta}}}_n^\tau)) 
= \widetilde{{\boldsymbol{\Delta}}^*}^T  \boldsymbol{A}_n^\gamma({\boldsymbol{\theta}}_0)\widetilde{{\boldsymbol{\Delta}}^*} 
+ 2 \widetilde{{\boldsymbol{\Delta}}^*}^T \boldsymbol{A}_\gamma({\boldsymbol{\theta}}_0) 
\sqrt{n}({\boldsymbol{\theta}}_n^\tau - {\boldsymbol{\theta}}_n^*) 
\nonumber\\
&-& 2 \widetilde{{\boldsymbol{\Delta}}^*}^T \boldsymbol{A}_\gamma({\boldsymbol{\theta}}_0) 
\sqrt{n}(\widetilde{{\boldsymbol{\theta}}}_n^\tau - \widetilde{{\boldsymbol{\theta}}_n^*})  
+ \sqrt{n}({\boldsymbol{\theta}}_n^\tau - {\boldsymbol{\theta}}_n^*)^T \boldsymbol{A}_\gamma({\boldsymbol{\theta}}_0) \sqrt{n}({\boldsymbol{\theta}}_n^\tau - {\boldsymbol{\theta}}_n^*)  
\nonumber\\
&-& 2 \sqrt{n}(\widetilde{{\boldsymbol{\theta}}}_n^\tau - \widetilde{{\boldsymbol{\theta}}_n^*})^T 
\boldsymbol{A}_\gamma({\boldsymbol{\theta}}_0)\sqrt{n}({\boldsymbol{\theta}}_n^\tau - {\boldsymbol{\theta}}_n^*)  
+ \sqrt{n}(\widetilde{{\boldsymbol{\theta}}}_n^\tau - \widetilde{{\boldsymbol{\theta}}_n^*})^T \boldsymbol{A}_\gamma({\boldsymbol{\theta}}_0) 
\sqrt{n}(\widetilde{{\boldsymbol{\theta}}}_n^\tau - \widetilde{{\boldsymbol{\theta}}_n^*}) 
+ o_P(1) +  o\left(1\right)\nonumber\\
&=& \boldsymbol{W}_n^T\boldsymbol{A}_\gamma({\boldsymbol{\theta}}_0)\boldsymbol{W}_n + o_P(1) +  o\left(1\right), \nonumber
\end{eqnarray}
where $\boldsymbol{W}_n = \left[\widetilde{{\boldsymbol{\Delta}}^*} + \sqrt{n}({\boldsymbol{\theta}}_n^\tau - {\boldsymbol{\theta}}_n^*) 
+\sqrt{n}(\widetilde{{\boldsymbol{\theta}}}_n^\tau - \widetilde{{\boldsymbol{\theta}}_n^*}) \right]$.
Thus  the asymptotic distribution of $S_\gamma({{\boldsymbol{\theta}}_n^\tau}, {\widetilde{{\boldsymbol{\theta}}}_n^\tau}) 
= 2 \sum_{i=1}^n d_{\gamma}(f_i(\cdot;{\boldsymbol{\theta}}_n^\tau), f_i(\cdot;{\widetilde{{\boldsymbol{\theta}}}_n^\tau})) $ 
under $\mathbf{\underline{F}}_{n,\epsilon,\mathbf{t}}^P$ is the same as the distribution of 
$\boldsymbol{W}^T \boldsymbol{A}_\gamma({\boldsymbol{\theta}}_0)\boldsymbol{W}$, 
where $\boldsymbol{W}$ is the asymptotic limit of $\boldsymbol{W}_n$.
But, from Results 1 and 2 of the Online Supplement one can show that, 
under $\mathbf{\underline{F}}_{n,\epsilon,\mathbf{t}}^P$, 
$\boldsymbol{W}\sim N_p(\widetilde{{\boldsymbol{\Delta}}^*},\widetilde{{\boldsymbol{\Sigma}}}_\tau({\boldsymbol{\theta}}_0))$.
This completes the proof of Part (i).

For Part (ii), consider the spectral decomposition of 
$\widetilde{{\boldsymbol{\Sigma}}}_\tau({\boldsymbol{\theta}}_0)^{1/2}
\boldsymbol{A}_\gamma({\boldsymbol{\theta}}_0)\widetilde{{\boldsymbol{\Sigma}}}_\tau({\boldsymbol{\theta}}_0)^{1/2}$
as
$
\widetilde{{\boldsymbol{\Sigma}}}_\tau({\boldsymbol{\theta}}_0)^{1/2}
\boldsymbol{A}_\gamma({\boldsymbol{\theta}}_0)\widetilde{{\boldsymbol{\Sigma}}}_\tau({\boldsymbol{\theta}}_0)^{1/2}
=\widetilde{\boldsymbol{P}}_{\tau,\gamma}(\boldsymbol{\theta}_0)^T 
\boldsymbol\Gamma_r \widetilde{\boldsymbol{P}}_{\tau,\gamma}(\boldsymbol{\theta}_0),
$
where $\widetilde{\boldsymbol{P}}_{\tau,\gamma}(\boldsymbol{\theta}_0)$ 
is as defined in 
the theorem and 
$\boldsymbol\Gamma_r=Diag\left\{\widetilde{\zeta_{i}^{\gamma,\tau}}(\boldsymbol\theta_0) : i=1,\ldots, r\right\}$.
Then $\boldsymbol{W}^T \boldsymbol{A}_\gamma({\boldsymbol{\theta}}_0)\boldsymbol{W}$ can be expressed as 
\begin{eqnarray}
&&\boldsymbol{W}^T\widetilde{{\boldsymbol{\Sigma}}}_\tau({\boldsymbol{\theta}}_0)^{-1/2} 
\left[\widetilde{{\boldsymbol{\Sigma}}}_\tau({\boldsymbol{\theta}}_0)^{1/2}\boldsymbol{A}_\gamma({\boldsymbol{\theta}}_0)
\widetilde{{\boldsymbol{\Sigma}}}_\tau({\boldsymbol{\theta}}_0)^{1/2}\right]
\widetilde{{\boldsymbol{\Sigma}}}_\tau({\boldsymbol{\theta}}_0)^{-1/2}\boldsymbol{W} \nonumber\\
&=&  \boldsymbol{W}^T\widetilde{{\boldsymbol{\Sigma}}}_\tau({\boldsymbol{\theta}}_0)^{-1/2} 
\left[\widetilde{\boldsymbol{P}}_{\tau,\gamma}(\boldsymbol{\theta}_0)^T 
\boldsymbol\Gamma_r \widetilde{\boldsymbol{P}}_{\tau,\gamma}(\boldsymbol{\theta}_0)\right]
\widetilde{{\boldsymbol{\Sigma}}}_\tau({\boldsymbol{\theta}}_0)^{-1/2}\boldsymbol{W} 
=(\boldsymbol{W}^\ast)^T\boldsymbol{\Gamma}_r\boldsymbol{W}^\ast, \nonumber
\end{eqnarray} 
where $\boldsymbol{W}^* = \widetilde{\boldsymbol{P}}_{\tau,\gamma}(\boldsymbol{\theta}_0)
\widetilde{{\boldsymbol{\Sigma}}}_\tau({\boldsymbol{\theta}}_0)^{-1/2}\boldsymbol{W}\sim 
N_p(\widetilde{\boldsymbol{\delta}}, \boldsymbol{I}_p)$
with $\widetilde{\boldsymbol{\delta}}=\left({\widetilde{\delta}_1}, \ldots, {\widetilde{\delta}_p}\right)^T$. 
This completes the proof of (ii).

Part (iii) follows from Part (i) using the series expansion of the distribution function of 
a linear combination of independent non-central chi-squares in terms of central chi-square distribution functions
as given in \cite{Kotz/etc:1967b}.
\hfill{$\square$}

%

\begin{corollary}\label{COR:8contiguous_power_composite}
Under the assumptions of Theorem \ref{THM:8asymp_Contaminated_power_composite},
we have the following.	
\begin{enumerate}
	\item ($\epsilon=0$): Asymptotic power under the contiguous alternatives $H_{1,n}$ is,	
	\begin{eqnarray}
		P_{\tau, \gamma}^*(\boldsymbol{\Delta},0; \alpha) 
		&=& \sum\limits_{v=0}^{\infty}~\widetilde{C_v^{\gamma,\tau}}(\boldsymbol{\theta}_0, {\boldsymbol{\Delta}})
		P\left(\chi_{r+2v}^2 > {s_\alpha^{\tau, \gamma}}/{\zeta_{(1)}^{\gamma,\tau}(\boldsymbol{\theta}_0)}\right).\nonumber
	\end{eqnarray}
	
	\item ($\boldsymbol{\Delta}=\boldsymbol{0}_p$): Asymptotic level under the contaminated distribution 
	$\mathbf{\underline{F}}_{n,\epsilon,\mathbf{t}}^L$ is 
	\begin{eqnarray}
	\alpha_\epsilon = P_{\tau, \gamma}^*(\boldsymbol{0}_p,\epsilon; \alpha)  
	= \sum\limits_{v=0}^{\infty}~\widetilde{C_v^{\gamma,\tau}}(\boldsymbol{\theta}_0, 
	\epsilon \boldsymbol{D}_\tau(\mathbf{t}, \boldsymbol{\theta}_0))
	P\left(\chi_{r+2v}^2 > {s_\alpha^{\tau, \gamma}}/{\zeta_{(1)}^{\gamma,\tau}(\boldsymbol{\theta}_0)}\right).
	\nonumber
	\end{eqnarray}
\end{enumerate}	
\end{corollary}



The following theorem then presents the PIF and LIF of the test in (\ref{EQ:8DPDTS2}).

\begin{theorem}\label{THM:8IF_power_composite}
Under the assumptions of Theorem \ref{THM:8asymp_Contaminated_power_composite}, 
if $\boldsymbol{D}_\tau(\mathbf{t}; {\theta_0})$ is bounded in $\mathbf{t}$,
then the power and level influence functions of the DPD based test (\ref{EQ:8DPDTS2}) are
	\begin{eqnarray}
		PIF(\mathbf{t}; S_{\gamma,\tau}, \mathbf{\underline{F}}_{\theta_0}) &=& 
		\boldsymbol{D}_\tau(\mathbf{t}, \boldsymbol{\theta}_0)^T 
		\widetilde{\boldsymbol{K}_{\gamma,\tau}}(\boldsymbol{\theta}_0,\boldsymbol{\Delta},\alpha),
		\mbox{and }~
		LIF(\mathbf{t}; S_{\gamma,\tau}, \mathbf{\underline{F}}_{\theta_0}) = 
		\boldsymbol{D}_\tau(\mathbf{t}, \boldsymbol{\theta}_0)^T 
		\widetilde{\boldsymbol{K}_{\gamma,\tau}}(\boldsymbol{\theta}_0,\boldsymbol{0}_p,\alpha),
		\nonumber
	\end{eqnarray} 
	where $
\widetilde{\boldsymbol{K}_{\gamma,\tau}}(\boldsymbol{\theta}_0,\boldsymbol{\Delta},\alpha) 
= \left(\sum\limits_{v=0}^{\infty}~ \left[\left.\frac{\partial}{\partial \boldsymbol{d}}
\widetilde{C_v^{\gamma,\tau}}(\boldsymbol{\theta}_0, \boldsymbol{d})\right|_{\boldsymbol{d}=\boldsymbol{\Delta}}\right]
P\left(\chi_{r+2v}^2 > {s_\alpha^{\tau, \gamma}}/{\zeta_{(1)}^{\gamma,\tau}(\boldsymbol{\theta}_0)}\right)\right).
	$
\end{theorem}
\noindent
\textbf{Proof}
Starting  with the  expression of $P_{\tau,\gamma}^*(\boldsymbol{\Delta}, \epsilon; \alpha)$ 
from Theorem \ref{THM:8asymp_Contaminated_power_composite},  we get
\begin{eqnarray}
&&PIF(\mathbf{t}; S_{\gamma,\tau}, \mathbf{\underline{F}}_{{\boldsymbol{\theta}}_0}) 
= \left.\frac{\partial}{\partial\epsilon}  P_{\tau,\gamma}^*({\boldsymbol{\Delta}}, \epsilon; \alpha)\right|_{\epsilon=0} 
=\sum\limits_{v=0}^{\infty}~\left.\frac{\partial}{\partial\epsilon}
\widetilde{C_v^{\gamma,\tau}}({\boldsymbol{\theta}}_0,\widetilde{{\boldsymbol{\Delta}}^*})\right|_{\epsilon=0} 
P\left(\chi_{r+2v}^2 > {s_\alpha^{\tau,\gamma}}/{\widetilde{\zeta_{(1)}^{\gamma,\tau}}({\boldsymbol{\theta}}_0)}
\right).
\label{EQ:7PIF_0simpleTest1}
\end{eqnarray} 
Now, for each $v \geq 0$, $\widetilde{C_v^{\gamma, \tau}}({\boldsymbol{\theta}}_0, \widetilde{{\boldsymbol{\Delta}}^*})$
depends on $\epsilon$ only through $\widetilde{{\boldsymbol{\Delta}}^*} 
= \left[ {\boldsymbol{\Delta}} + \epsilon \boldsymbol{D}_\tau(\mathbf{t}, {{\boldsymbol{\theta}}_0})\right]$.
Consider a Taylor series expansion of $\widetilde{C_v^{\gamma, \tau}}({\boldsymbol{\theta}}_0, \boldsymbol{d})$ 
with respect to $\boldsymbol{d}$ around $\boldsymbol{d}={\boldsymbol{\Delta}}$
and evaluate it at $\boldsymbol{d}=\widetilde{{\boldsymbol{\Delta}}^*}$ to get
\begin{eqnarray}
&& \widehat{C_v^{\gamma, \tau}}({\boldsymbol{\theta}}_0, \widetilde{{\boldsymbol{\Delta}}^*}) 
= \widetilde{C_v^{\gamma, \tau}}({\boldsymbol{\theta}}_0, {{\boldsymbol{\Delta}}}) 
+ (\widetilde{{\boldsymbol{\Delta}}^*} - {\boldsymbol{\Delta}})^T  \left[\left.\frac{\partial}{\partial \boldsymbol{d}}
\widetilde{C_v^{\gamma,\tau}}({\boldsymbol{\theta}}_0,\boldsymbol{d})^T\right|_{\boldsymbol{d}={\boldsymbol{\Delta}}}\right]
+ o(||\widetilde{{\boldsymbol{\Delta}}^*} - {\boldsymbol{\Delta}}||)\nonumber\\
&&= \widetilde{C_v^{\gamma, \tau}}({\boldsymbol{\theta}}_0, {{\boldsymbol{\Delta}}}) 
+ \epsilon \boldsymbol{D}(\mathbf{t}, {{\boldsymbol{\theta}}_0})^T \cdot \left[\left.\frac{\partial}{\partial \boldsymbol{d}}
\widetilde{C_v^{\gamma,\tau}}({\boldsymbol{\theta}}_0,\boldsymbol{d})\right|_{\boldsymbol{d}={\boldsymbol{\Delta}}}\right]
+ o(\epsilon ||\boldsymbol{D}(\mathbf{t}, {{\boldsymbol{\theta}}_0})||).\nonumber
\end{eqnarray}
Now, since $\boldsymbol{D}(\mathbf{t}, {{\boldsymbol{\theta}}_0})$ is finite, 
differentiating it with respect to $\epsilon$ and evaluating at $\epsilon=0$, 
we get that
$
\left.\frac{\partial}{\partial\epsilon}
\widetilde{C_v^{\gamma,\tau}}({\boldsymbol{\theta}}_0,\widetilde{{\boldsymbol{\Delta}}^*})\right|_{\epsilon=0}$
$=$$ \boldsymbol{D}(\mathbf{t}, {{\boldsymbol{\theta}}_0})^T \left[\left.\frac{\partial}{\partial \boldsymbol{d}}
\boldsymbol{C_v^{\gamma,\tau}}({\boldsymbol{\theta}}_0,\boldsymbol{d})\right|_{\boldsymbol{d}={\boldsymbol{\Delta}}}\right].
$
Combining it with Equation (\ref{EQ:7PIF_0simpleTest1}), we finally get the required PIF.
The LIF is then obtained from the PIF by substituting $\boldsymbol{\Delta}=\boldsymbol{0}_p$. 
\hfill{$\square$}

Note that, under the general INH set-up, both LIF and PIF are bounded 
whenever the IFs of the MDPDE under the null and overall parameter space are bounded. 
But this is the case for most statistical models at $\tau>0$
implying the size and power robustness of the corresponding DPD based tests.

\section{Application: Testing General Linear Hypothesis under the Normal Linear Regression}
\label{SEC:linear_regression}

\noindent
We assume that, given fixed covariates ${\boldsymbol{x}_1}, \ldots, {\boldsymbol{x}_n}\in \mathbb{R}^p$,  
the (random) responses $y_1, \ldots, y_n$ satisfy the relation
\begin{equation}
y_i = {\boldsymbol{x}_i}^T{{\mathbf {\boldsymbol{\beta}}}} + \epsilon_i, ~~i = 1, \ldots, n,
\end{equation}
where $\epsilon_i$'s are independent and identically distributed as $N(0, \sigma^2)$
and ${{\mathbf {\boldsymbol{\beta}}}}= ({{\beta}}_1, \ldots, {{\beta}}_p)^T$ is the vector of regression coefficients.
Thus, $y_i$s are INH with   $y_i \sim N({\boldsymbol{x}_i}^T{{\mathbf {\boldsymbol{\beta}}}}, \sigma^2)$ for each $i$. 
The most common general linear hypothesis is given by 
\begin{eqnarray}
	H_0 : \boldsymbol{L}^T\boldsymbol{\beta} = \boldsymbol{l}_0 ~~~~\mbox{ against }~~~~~ 
	H_1 : \boldsymbol{L}^T\boldsymbol{\beta} \neq \boldsymbol{l}_0, ~~
	\label{EQ:9Gen_lin_hypothesis}
\end{eqnarray}  
where $\sigma$  is unknown in both cases, 
$\boldsymbol{L}$ is a $p\times r$ known matrix ($r\leq p$) and $\boldsymbol{l}_0$ is a known $p$-vector of reals. 
We assume that $rank(\boldsymbol{L})=r$ so that the null hypothesis in (\ref{EQ:9Gen_lin_hypothesis})
is feasible with solution $\boldsymbol{\beta}_0$ and also of the form (\ref{EQ:8composite_hypo})
with $\boldsymbol{\theta} \in \Theta_0 = \left\{\boldsymbol\beta_0 \in \mathbb{R}^p : \boldsymbol{L}^T\boldsymbol\beta_0=\boldsymbol{l}_0\right\}
\times [0, \infty)\subset \Theta=\mathbb{R}^p\times [0, \infty)$.

To define the DPD based test for testing (\ref{EQ:9Gen_lin_hypothesis}), let 
$\widetilde{\boldsymbol{\theta}}_n^\tau = (\widetilde{\boldsymbol{\beta}}_n^\tau, \widetilde{\sigma}_n^\tau)$ 
and $\boldsymbol{\theta}_n^\tau = (\boldsymbol{\beta}_n^\tau, \sigma_n^\tau)$
denote the RMDPDE of $\boldsymbol{\theta} = (\boldsymbol{\beta}, \sigma)$ under $H_0$ in (\ref{EQ:9Gen_lin_hypothesis})
and their unrestricted MDPDE, respectively, both with tuning parameter $\tau$. 
Note that, $\widetilde{\boldsymbol{\beta}}_n^\tau = \boldsymbol{\beta}_0$ and 
hence our DPD based test statistics (\ref{EQ:8DPDTS2}) for testing (\ref{EQ:9Gen_lin_hypothesis}) becomes
\begin{eqnarray}
	S_{\gamma}( \boldsymbol{\theta}_n^\tau, \widetilde{\boldsymbol{\theta}}_n^\tau)  =
	\frac{2\sqrt{1+\gamma}}{\gamma(\sqrt{2\pi}\widetilde{\sigma}_n^\tau)^\gamma}
	\left[n C_1 - C_2\sum_{i=1}^n ~ 
	e^{-\frac{\gamma(\boldsymbol{\beta}_n^\tau-\boldsymbol{\beta}_0)^T(\boldsymbol{x}_i\boldsymbol{x}_i^T)
			(\boldsymbol{\beta}_n^\tau-\boldsymbol{\beta}_0)}{2(\gamma(\sigma_n^\tau)^2 + 
			(\widetilde{\sigma}_n^\tau)^2)}} \right],\mbox{ for }\gamma>0,\nonumber
	\label{EQ:9composite_Test_stat_normal}\\
	S_{0}(  \boldsymbol{\theta}_n^\tau, \widetilde{\boldsymbol{\theta}}_n^\tau)  = n
	\left[\log\left(\frac{(\widetilde{\sigma}_n^\tau)^2}{(\sigma_n^\tau)^2}\right) - 1
	+ \frac{(\sigma_n^\tau)^2}{(\widetilde{\sigma}_n^\tau)^2} \right] 
	+ \frac{(\boldsymbol{\beta}_n^{\tau}-\boldsymbol{\beta}_0)^T
		(\boldsymbol{X}^T\boldsymbol{X})(\boldsymbol{\beta}_n^\tau-\boldsymbol{\beta}_0)}{(\widetilde{\sigma}_n^\tau)^2},
	\nonumber\label{EQ:9composite_Test_stat0_normal}
\end{eqnarray}
with $C_1 = [\gamma(\sigma_n^\tau)^\gamma + (\widetilde{\sigma}_n^\tau)^\gamma]
(1+\gamma)^{-1}    (\sigma_n^\tau)^{-\gamma}$,
$C_2 =  \sigma_n^\tau \sqrt{1+\gamma}[\gamma(\sigma_n^\tau)^2 + (\widetilde{\sigma}_n^\tau)^2]^{-1/2}$ 
and $\boldsymbol{X}^T = \left[\boldsymbol{x}_1, \ldots, \boldsymbol{x}_n\right]_{p\times n}$.
At $\gamma=\tau=0$, it coincides with the LRT statistic. 

In the following, we derive the properties of this DPD based test under the general linear hypothesis (\ref{EQ:9Gen_lin_hypothesis}); later we
and illustrate their applications for the example of testing for the first $r\leq p$ components of $\boldsymbol{\beta}$.

%
\noindent
\textit{Asymptotic Distributions:}

The asymptotic distribution of the MDPDE $\boldsymbol{\theta}_n^\tau = (\boldsymbol{\beta}_n^\tau, \sigma_n^\tau)$ 
under this fixed-design linear regression model (LRM) has been derived in \cite{Ghosh/Basu:2013}; 
under Assumptions (R1)--(R2) of the Online Supplement, 
if $\boldsymbol{\theta}_0=(\boldsymbol{\beta}_0, \sigma_0)$ is the true parameter value, 
the the MDPDEs $\widehat{{{\mathbf {\boldsymbol{\beta}}}}}$ and $\widehat\sigma^2$ are both consistent 
and asymptotically independent with 
$(\boldsymbol{X}^T\boldsymbol{X})^{\frac{1}{2}}(\boldsymbol{\beta}_n^\tau- {\boldsymbol{\beta}}_0)
\displaystyle\mathop{\rightarrow}^\mathcal{D} 
N_p\left(\boldsymbol{0}, \upsilon_\tau^{{\mathbf {\boldsymbol{\beta}}}}\boldsymbol{I}_p\right)$
and $\sqrt{n}(\widehat\sigma^2 - \sigma_0^2) \displaystyle\mathop{\rightarrow}^\mathcal{D} N(0, \upsilon_\tau^e)$, 
where 
$\upsilon_\tau^{{\mathbf {\boldsymbol{\beta}}}}= \sigma_0^2 \left(1 + \frac{\tau^2}{1+2\tau}\right)^{\frac{3}{2}}$
and
$\upsilon_\tau^e = \frac{4\sigma_0^4}{(2+\tau^2)^2} 
\left[2(1+2\tau^2)\left(1+ \frac{\tau^2}{1+2\tau}\right)^{\frac{5}{2}} - \tau^2(1+\tau)^2\right]$.

The asymptotic distribution of the RMDPDE  
$\widetilde{\boldsymbol{\theta}}_n^\tau =(\widetilde{\boldsymbol{\beta}}_n^\tau, \widetilde{\sigma}_n^\tau)$
can be obtained from Result 2 of the Online Supplement with 
$\boldsymbol{\upsilon}(\boldsymbol{\beta}, \sigma) = \boldsymbol{L}^T\boldsymbol{\beta} - \boldsymbol{\beta}_0$, 
$\boldsymbol{\Upsilon}(\boldsymbol{\beta}, \sigma) = [\boldsymbol{L}^T~\boldsymbol{0}_r]^T
$ and $\nabla^2H_n(\boldsymbol{\theta}) = (1+\tau) \boldsymbol{A}_n^\tau(\boldsymbol{\theta})$,
which is presented in the following Theorem. 
Note that Assumptions (R1)--(R2) imply Assumptions (A1)--(A7) under any $\boldsymbol{\theta} \in \Theta$ in the LRM
and hence for $\boldsymbol{\theta}\in \Theta_0$ \citep[Lemma 6.1]{Ghosh/Basu:2013}.

\begin{theorem}
Suppose $rank(\boldsymbol{L})=r$, Assumptions (R1)--(R2) of the Online Supplement hold and 
the true parameter value $(\boldsymbol{\beta}_0, \sigma_0)\in \Theta_0$.
Then, for $\tau \geq 0$, 
there exists consistent RMDPDE
$(\widetilde{\boldsymbol{\beta}}_n^\tau, \widetilde{\sigma}_n^\tau)$ under $H_0$ in (\ref{EQ:9Gen_lin_hypothesis})
which are asymptotically independent and
$(\boldsymbol{X}^T \boldsymbol{X})^{\frac{1}{2}}\widetilde{\boldsymbol{P}_n}^{-1}
(\widetilde{\boldsymbol{\beta}}_n^\tau- {\boldsymbol{\beta}_0})
\displaystyle\mathop{\rightarrow}^\mathcal{D}  
N_p\left(\boldsymbol{0}_p, \upsilon_\tau^{\boldsymbol{\beta}} \boldsymbol{I}_p\right)$ 
and 
$\sqrt{n}\left[(\widetilde{\sigma}_n^\tau)^2 - \sigma_0^2\right]
\displaystyle\mathop{\rightarrow}^\mathcal{D}$ $N(0, \upsilon_\tau^e)$,
where
$
\widetilde{\boldsymbol{P}_n} = \left[\boldsymbol{I}_p - \boldsymbol{L}\{\boldsymbol{L}^T(\boldsymbol{X}^T\boldsymbol{X})^{-1}\boldsymbol{L}\}^{-1}\boldsymbol{L}^T(\boldsymbol{X}^T\boldsymbol{X})^{-1}\right].
$ 
%
%
\label{THM:9asymp_RMDPDE_lin_reg}
\end{theorem}

Note that, the asymptotic relative efficiency of the RMDPDEs of $\boldsymbol{\beta}$ and $\sigma^2$ are exactly
the same as that of their unrestricted versions, which are quite high for small $\tau>0$ \citep{Ghosh/Basu:2013}. 

Our next theorem presents the asymptotic null distribution of the DPD test statistics in the LRM; 
its proof follows from Result 3 of the Online Supplement.

\begin{theorem}
Suppose $rank(\boldsymbol{L})=r$, Assumptions (R1)--(R3) of the Online Supplement  hold and 
the true parameter value $(\boldsymbol{\beta}_0, \sigma_0)\in \Theta_0$.
Then, the asymptotic distribution of $S_{\gamma}( \boldsymbol{\theta}_n^\tau, \widetilde{\boldsymbol{\theta}}_n^\tau)$
under $H_0$ in (\ref{EQ:9Gen_lin_hypothesis})  coincides with the distribution of 
$\zeta_1^{\gamma, \tau} \sum_{i=1}^r  \lambda_i Z_i^2,$
where $Z_1, \cdots,Z_r$ are independent standard normal variables, $\lambda_1, \cdots, \lambda_r$ are nonzero eigenvalues of 
$\boldsymbol{Q}_x = \left(\boldsymbol{L}
\left[\boldsymbol{L}^T\boldsymbol\Sigma_x^{-1}\boldsymbol{L}\right]^{-1}\boldsymbol{L}^T\boldsymbol\Sigma_x^{-1}\right)$
and  $\zeta_1^{\gamma, \tau}=(1+\gamma) s_\gamma \upsilon_\tau^{{\mathbf {\boldsymbol{\beta}}}}$ with 
$s_\gamma= (2\pi)^{-\frac{\gamma}{2}} \sigma^{-(\gamma+2)}(1+\gamma)^{-\frac{3}{2}}$.
	\label{THM:9asymp_null_composite_linReg}
\end{theorem}

Further, from the general theory from \cite{Ghosh/Basu:2017}, 
this DPD based test is consistent at any fixed alternative.
Under the assumptions of Theorem \ref{THM:9asymp_null_composite_linReg}, 
the asymptotic distribution of  $S_{\gamma}( {{\boldsymbol{\theta}}_n^\tau}, {\widetilde{{\boldsymbol{\theta}}}_n^\tau})$ 
under $H_{1,n}': {\boldsymbol{\beta}}={\boldsymbol{\beta}}_n$, ${\boldsymbol{\beta}}_n = {\boldsymbol{\beta}}_0 + n^{-\frac{1}{2}}\boldsymbol\Delta_1$,
is the same as that of $\zeta_1^{\gamma, \tau} \sum_{i=1}^r  \lambda_i (Z_i +\delta_i)^2,$ 
where
$
\left({\delta_1}, \cdots, {\delta_p}\right)^T
= \widetilde{\boldsymbol{N}}
\left[\upsilon_\tau^{\boldsymbol{\beta}}\boldsymbol\Sigma_x^{-1}\boldsymbol{Q}_x
\right]^{-1/2}\boldsymbol\Delta_1,
$
with $\widetilde{\boldsymbol{N}}$ being the matrix of normalized eigenvectors of $\boldsymbol{Q}_x$
(Theorem \ref{THM:8asymp_Contaminated_power_composite} at $\epsilon=0$).  
This leads to the asymptotic contiguous power which decreases as $\tau = \gamma$ increases.

\noindent
\textit{Influence Functions:}\\
From Section \ref{SEC:8composite_testing},
the first order IF of the DPD based test is always zero when evaluated at $H_0$
and its second order IF, given by (\ref{EQ:IF2_TS}), depends  on the IFs of the MDPDE functional, 
say $\boldsymbol{U}_\tau=(\boldsymbol{U}_\tau^{\boldsymbol{\beta}}, \boldsymbol{U}_\tau^{\sigma})$, 
and the RMDPDE functional, say 
$\widetilde{\boldsymbol{U}}_\tau = (\widetilde{\boldsymbol{U}}_\tau^{\boldsymbol{\beta}}, 
\widetilde{\boldsymbol{U}}_\tau^{\sigma})$, of $\boldsymbol{\theta}=(\boldsymbol{\beta}, \boldsymbol{\sigma})$.
The IF of  $\boldsymbol{U}_\tau$ has already been derived in \cite{Ghosh/Basu:2013}. 
Under contamination in all directions, 
the IFs of $\boldsymbol{U}_\tau^{\boldsymbol{\beta}}$ and $\boldsymbol{U}_\tau^{\sigma}$,
at $\mathbf{\underline{G}} = \mathbf{\underline{F}}_{\boldsymbol{\theta}_0}$, are individually given by
\begin{eqnarray}
\mathcal{IF}(\boldsymbol{t}, \boldsymbol{U}_\tau^{{{\boldsymbol{\beta}}}}, \mathbf{\underline{F}}_{\boldsymbol{\theta}_0}) &=&
(1+\tau)^{\frac{3}{2}} (\boldsymbol{X}^T\boldsymbol{X})^{-1}
\sum_{i=1}^n {\boldsymbol{x}_{i}} (t_{i} - {\boldsymbol{x}_{i}}^T{{\mathbf {\boldsymbol{\beta}}}})
e^{-\frac{\tau(t_{i} - {\boldsymbol{x}_{i}}^T{{\mathbf {\boldsymbol{\beta}}}})^2}{2\sigma^2}},\nonumber
\label{EQ:IF_i0_reg}\\
%
\mathcal{IF}(\boldsymbol{t}, \boldsymbol{U}_\tau^\sigma, \mathbf{\underline{F}}_{\boldsymbol{\theta}_0}) &=&
\frac{2(1+\tau)^{\frac{5}{2}}}{n(2+\tau^2)}\sum_{i=1}^n \left\{ (t_{i} -{\boldsymbol{x}_{i}}^T{\boldsymbol{\beta}})^2 
- \sigma^2 \right\}e^{-\frac{\tau(t_{i} - {\boldsymbol{x}_{i}}^T{\boldsymbol{\beta}})^2}{2\sigma^2}} 
+ \frac{2\tau(1+\tau)^2}{(2+\tau^2)}.\nonumber
\label{EQ:IF_i0_var}
\end{eqnarray}
Now we derive the IF of the RMDPDE $\widetilde{\boldsymbol{U}}_\tau = (\widetilde{\boldsymbol{U}}_\tau^{\boldsymbol{\beta}}, 
\widetilde{\boldsymbol{U}}_\tau^{\sigma})$ following the general theory of \cite{Ghosh/Basu:2017}.
It follows that, under contamination in all directions, the IFs of $\widetilde{\boldsymbol{U}}_\tau^{\boldsymbol{\beta}}$
and $\widetilde{\boldsymbol{U}}_\tau^\sigma$ are also independently obtainable at 
$\mathbf{\underline{G}} = \mathbf{\underline{F}}_{\theta_0}$ as
$\mathcal{IF}(\boldsymbol{t}, \widetilde{\boldsymbol{U}}_\tau^\sigma, \mathbf{\underline{F}}_{\boldsymbol{\theta}_0})
=\mathcal{IF}(\boldsymbol{t}, \boldsymbol{U}_\tau^\sigma, \mathbf{\underline{F}}_{\boldsymbol{\theta}_0})
$
and  
\begin{eqnarray}
\mathcal{IF}(\boldsymbol{t}, \widetilde{\boldsymbol{U}}_\tau^{\boldsymbol{\beta}}, 
{\underline{\boldsymbol{F}}}_{\boldsymbol{\theta}_0}) &=&
	\left[\boldsymbol{\Psi}_{1,n}^{\tau, 0}(\boldsymbol{\beta})^T\boldsymbol{\Psi}_{1,n}^{\tau, 0}(\boldsymbol{\beta}) 
	+ \boldsymbol{L}\boldsymbol{L}^T\right]^{-1}\boldsymbol{\Psi}_{1,n}^{\tau, 0}(\boldsymbol{\beta})^T
\frac{1}{n}\sum_{i=1}^n\left\{\boldsymbol{u}_i^{(0)}(t_i,\boldsymbol{\beta}) 
	\phi(t_i; \boldsymbol{x}_i^T\boldsymbol{\beta}, \sigma)^{\tau}- \boldsymbol{\xi}_i^{(0)}(\boldsymbol{\beta}_0)\right\},
	\label{EQ:9IF_i0_RMDPDE_reg}\nonumber
\end{eqnarray}
where $\phi(y; \mu, \sigma)$ denotes the density of $N(\mu, \sigma^2)$ at $y$, 
$\boldsymbol{\xi}_i^{(0)}(\boldsymbol{\beta}) = \int \boldsymbol{u}_i^{(0)}(y,\boldsymbol{\beta})
\phi(y; \boldsymbol{x}_i^T\boldsymbol{\beta}, \sigma)^{1+\tau}dy$   and 
$\boldsymbol{\Psi}_{1,n}^{\tau, 0}(\boldsymbol{\beta}) = \frac{1}{n}\sum_{i=1}^n
\int \boldsymbol{u}_i^{(0)}(y,\boldsymbol{\beta})\boldsymbol{u}_i^{(0)}(y,\boldsymbol{\beta})^T
\phi(y; \boldsymbol{x}_i^T\boldsymbol{\beta}, \sigma)^{1+\tau}dy$
with $\boldsymbol{u}_i^{(0)}(y,\boldsymbol{\beta})$ being the likelihood score function of $\boldsymbol{\beta}$ 
under the restriction  of $H_0$ in (\ref{EQ:9Gen_lin_hypothesis}). 
Since the IF of error variance $\sigma^2$ under restrictions is the same 
as that in the unrestricted case, it follows from (\ref{EQ:IF2_TS}) that
the second order IF of the DPD based test statistic is
\begin{eqnarray}
\mathcal{IF}_2(\mathbf{t}, S_{\gamma, \tau}, \underline{\mathbf{F}}_{\theta_0}) 
= (1+\gamma)\zeta_\gamma  \boldsymbol{D}_{\tau}^{\boldsymbol{\beta}}(\mathbf{t}, \boldsymbol{\theta}_0)^T 
(\boldsymbol{X}^T\boldsymbol{X})\boldsymbol{D}_{\tau}^{\boldsymbol{\beta}}(\mathbf{t}, \boldsymbol{\theta}_0),\nonumber
\end{eqnarray}
with $\boldsymbol{D}_{\tau}^{\boldsymbol{\beta}}(\mathbf{t}, \boldsymbol{\theta}_0) =\left[
\mathcal{IF}(\boldsymbol{t}, {\boldsymbol{U}}_\tau^{\boldsymbol{\beta}}, {\underline{\boldsymbol{F}}}_{\boldsymbol{\theta}_0})
- \mathcal{IF}(\boldsymbol{t}, \widetilde{\boldsymbol{U}}_\tau^{\boldsymbol{\beta}}, 
{\underline{\boldsymbol{F}}}_{\boldsymbol{\theta}_0}) \right]$.
At $\tau>0$, this second order IF is bounded in $\boldsymbol{t}$ implying robustness.
The case $\tau=0$ is not conclusive; an example is provided later.

\noindent
\textit{Power and Level Robustness:}\\
It follows from Theorem \ref{THM:8IF_power_composite} that the asymptotic distribution of 
$S_{\gamma}( {\theta_n^\tau}, {\widetilde{\theta}_n^\tau})$  under $H_{1,n}'$ along with contiguous contamination
is given by  $\zeta_1^{\gamma, \tau} \sum_{i=1}^r  \lambda_i (Z_i +\widetilde{\delta}_i)^2,$ 
where
$
\left(\widetilde{\delta_1}, \cdots, \widetilde{\delta_p}\right)^T
= \widetilde{\boldsymbol{N}}
\left[\upsilon_\tau^{\boldsymbol{\beta}}\boldsymbol\Sigma_x^{-1}\boldsymbol{Q}_x
\right]^{-1/2}
\left[\boldsymbol{\Delta} + \epsilon \boldsymbol{D}_\tau^{\boldsymbol{\beta}}(\mathbf{t}, \boldsymbol{\theta}_0)\right],
$
under the assumptions of Theorem \ref{THM:9asymp_null_composite_linReg}.
Then, the PIF and LIF can be derived empirically from the infinite sum representation 
given in Theorem \ref{THM:8IF_power_composite}.   
However, for any  general restriction, 
both the LIF and PIF depend on the contamination points $\mathbf{t}$ 
only through the quantity $\boldsymbol{D}_\tau^{\boldsymbol{\beta}}(\mathbf{t}, \boldsymbol{\theta}_0)$,
which is independent of the IF of the estimates of $\sigma$
and hence independent of its robustness properties.
\noindent
\textit{Example \ref{SEC:linear_regression}.1 [Test for only the first $r\leq p$ components of $\beta$]:}\\
Let us now illustrate the above results for the most common case of (\ref{EQ:9Gen_lin_hypothesis}),
where we fix the first $r$ components ($r\leq p$) of $\boldsymbol{\beta}$ at a pre-fixed values $\boldsymbol{\beta}_0^{(1)}$. 
So, our null hypothesis becomes  $H_0: \boldsymbol{\beta}^{(1)}=\boldsymbol{\beta}_0^{(1)}$,
where $\boldsymbol{\beta}^{(1)}$ denote the first $r$-components of 
$\boldsymbol{\beta}=(\boldsymbol{\beta}^{(1)T}, \boldsymbol{\beta}^{(2)T})^T$.
%
In terms (\ref{EQ:9Gen_lin_hypothesis}), we have 
$\boldsymbol{L}=[\boldsymbol{I}_r ~ \boldsymbol{O}_{r\times (p-r)}]^T$
and $\boldsymbol{l}_0 = \boldsymbol{\beta}_0^{(1)}$.
Let us consider the partitions $\boldsymbol{\beta}_0^T=(\boldsymbol{\beta}_0^{(1)T}, \boldsymbol{\beta}_0^{(2)T})$, 
$\boldsymbol{x}_i^T=(\boldsymbol{x}_i^{(1)T}, \boldsymbol{x}_i^{(2)T})$ 
and $\boldsymbol{X}=[\boldsymbol{X}_1 ~ \boldsymbol{X}_2]$, 
where $\boldsymbol{\beta}_0^{(1)}$ and $\boldsymbol{x}_i^{(1)}$ are $r$-vectors 
and $\boldsymbol{X}_1$ is the $n\times r$ matrix consisting of the first $r$ columns of $\boldsymbol{X}$. 
Then, the distribution of the RMDPDEs of first $r$ fixed components of $\boldsymbol{\beta}$ 
turns out to be degenerate at their given values $\boldsymbol{\beta}_0^{(1)}$. 
We can derive the asymptotic distribution for rest of the components
using Theorem \ref{THM:9asymp_RMDPDE_lin_reg}, as given by
$(\boldsymbol{X}^T\boldsymbol{X})_{22.1}^{\frac{1}{2}}
[(\widetilde{\boldsymbol{\beta}}_n^\tau)^{(2)}- {\boldsymbol{\beta}}^{(2)}]
\displaystyle\mathop{\rightarrow}^\mathcal{D} 
N_{p-r}(\boldsymbol{0}_{p-r}, \upsilon_\tau^{\boldsymbol{\beta}}\boldsymbol{I}_{p-r})$,
where $(\boldsymbol{X}^T\boldsymbol{X})_{22.1} = [(\boldsymbol{X}_2^T\boldsymbol{X}_2) 
- (\boldsymbol{X}_2^T\boldsymbol{X}_1)(\boldsymbol{X}_1^T\boldsymbol{X}_1)^{-1}(\boldsymbol{X}_1^T\boldsymbol{X}_2)]$.
Next, considering the DPD based test for this problem, under the assumptions of Theorem \ref{THM:9asymp_null_composite_linReg},
the asymptotic null distribution of 
$S_{\gamma}( {\boldsymbol{\theta}_n^\tau}, {\widetilde{\boldsymbol{\theta}}_n^\tau})/\zeta_1^{\gamma, \tau}$ 
is simply chi-square with df $r$. So, the critical values are straightforward. 
%
%
In terms of robustness, the IF of the RMDPDE and the second order IF of the DPD based test statistics 
further simplify in this case as
\begin{equation}
\mathcal{IF}(\boldsymbol{t}, \widetilde{\boldsymbol{U}}_\tau^{\boldsymbol{\beta}}, \underline{\mathbf{F}}_{\theta_0}) =
\begin{bmatrix}
~~\boldsymbol{0}_r~\\
~~(1+\tau)^{\frac{3}{2}} (\boldsymbol{X}_2^T \boldsymbol{X}_2)^{-1}\sum_{i=1}^n \boldsymbol{x}_{i}^{(2)}
(t_{i} - \boldsymbol{x}_{i}^T{\boldsymbol{\beta}})
e^{-\frac{\tau(t_{i} - \boldsymbol{x}_{i}^T{\boldsymbol{\beta}})^2}{2\sigma^2}}~
\end{bmatrix},\nonumber
\label{EQ:9IF_i0_RMDPDE_reg_eg2}
\end{equation}
\begin{eqnarray}
\mathcal{IF}_2(\mathbf{t}, S_{\gamma, \tau}, \underline{\mathbf{F}}_{\theta_0}) 
= (1+\gamma)\zeta_\gamma (1+\tau)^{\frac{3}{2}} \sum_{i=1}^n 
\left[\boldsymbol{x}_{i}^{(1)T} \boldsymbol{M}_x \boldsymbol{x}_{i}^{(1)} \right] 
(t_{i} - \boldsymbol{x}_{i}^T{\boldsymbol{\beta}})^2 
e^{-\frac{\tau(t_{i} - \boldsymbol{x}_{i}^T{\boldsymbol{\beta}})^2}{\sigma^2}}, 
\nonumber
\end{eqnarray}
where 
$\boldsymbol{M}_x = (\boldsymbol{X}^T \boldsymbol{X})_{11.2}^{-1}(\boldsymbol{X}_1^T\boldsymbol{X}_1)
(\boldsymbol{X}^T \boldsymbol{X})_{11.2}^{-1},$
with $(\boldsymbol{X}^T\boldsymbol{X})_{11.2} = [(\boldsymbol{X}_1^T\boldsymbol{X}_1) 
- (\boldsymbol{X}_1^T\boldsymbol{X}_2)(\boldsymbol{X}_2^T\boldsymbol{X}_2)^{-1}(\boldsymbol{X}_2^T\boldsymbol{X}_1)]$.
In order to obtain the PIF, we consider the contiguous alternatives
$H_{1,n}'': \boldsymbol{\beta}^{(1)}=\boldsymbol{\beta}_n^{(1)}$, where $\boldsymbol{\beta}_n^{(1)} = \boldsymbol{\beta}_0^{(1)} 
+ \frac{\boldsymbol{\Delta}_1^{(1)}}{\sqrt{n}}$ and $\boldsymbol{\Delta}_1^{(1)}$ is the first $r$ components of 
$\boldsymbol{\Delta}_1=(\boldsymbol{\Delta}_1^{(1)}, \boldsymbol{\Delta}_1^{(2)})$.
Then, following Theorem \ref{THM:8IF_power_composite}, we get 
\begin{eqnarray}
&& PIF(\mathbf{t}; S_{\gamma, \tau}^{(1)}, \mathbf{\underline{F}}_{\theta_0})\label{EQ:9PIF_test_unknownVar}
= \widetilde{K_{\gamma,\tau}}
\left(\boldsymbol{\Delta}_1^{(1)T}\boldsymbol{\Sigma}_x^{(11)}\boldsymbol{\Delta}_1^{(1)}, r\right)
\sum\limits_{i=1}^{n}[\boldsymbol{\Delta}_1^{(1)T}\boldsymbol{x}_{i}^{(1)}] (t_{i} - \boldsymbol{x}_{i}^T{\boldsymbol{\beta}_0})
e^{-\frac{\tau(t_{i} - \boldsymbol{x}_{i}^T{\boldsymbol{\beta}_0})^2}{2\sigma_0^2}}.
\end{eqnarray} 
where $\boldsymbol{\Sigma}_x ^{(11)}$ is the $r\times r$ principle minor of $\boldsymbol{\Sigma}_x$.
Note that, as we have fixed the first $r$ components of $\beta$, their IFs are zero. 
Further, all these IFs are bounded whenever $\tau >0$ and unbounded at $\tau=0$. 
Thus the DPD based test with $\tau>0$ is stable in its asymptotic power but the LRT ($\tau=0$) is not. 

Finally, substituting $\Delta_1^{(1)} = 0$ in (\ref{EQ:9PIF_test_unknownVar}), 
we get $LIF(\mathbf{t}; S_{\gamma, \tau}^{(1)}, \mathbf{\underline{F}}_{\theta_0})=0$
for all $\tau  > 0$ implying robustness in terms of asymptotic level of the DPD based tests.


\begin{figure}[h]
	\centering
	\subfloat[$n=30$, $\tau=\gamma=0$]{
		\includegraphics[width=0.3\textwidth]{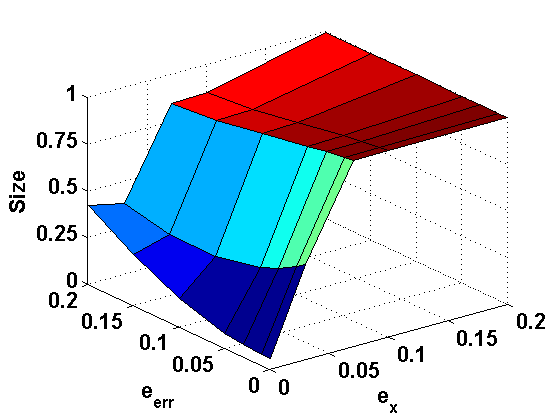}
		\label{FIG:9size_unknown_30_0}}
	~ 
	\subfloat[$n=30$, $\tau=\gamma=0.5$]{
		\includegraphics[width=0.3\textwidth]{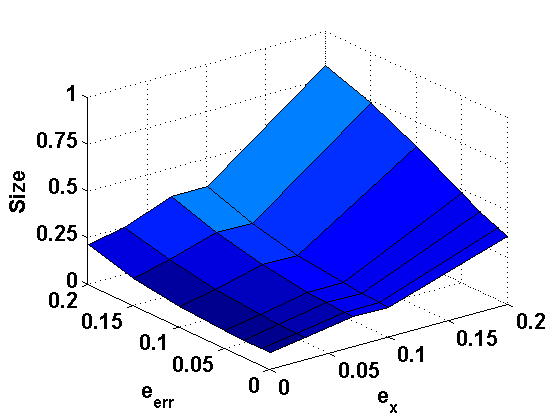}
		\label{FIG:9size_unknown_30_05}}
	~ 
	\subfloat[$n=30$, $\tau=\gamma=1$]{
		\includegraphics[width=0.3\textwidth]{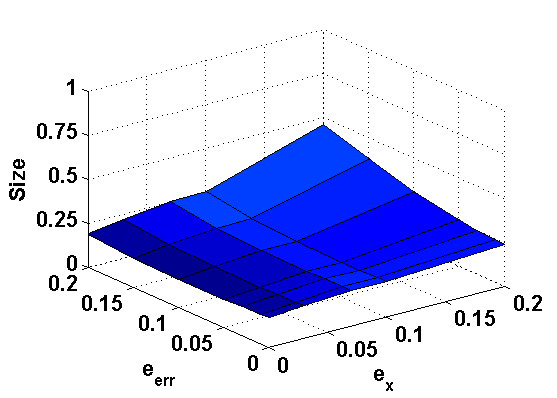}
		\label{FIG:9size_unknown_30_1}}
	\\
	\subfloat[$n=50$, $\tau=\gamma=0$]{
		\includegraphics[width=0.3\textwidth]{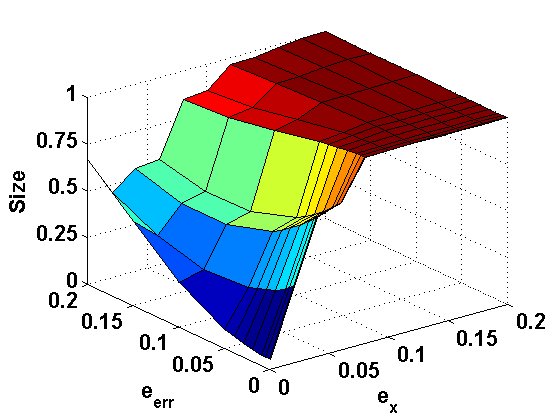}
		\label{FIG:9size_unknown_50_0}}
	~ 
	\subfloat[$n=50$, $\tau=\gamma=0.5$]{
		\includegraphics[width=0.3\textwidth]{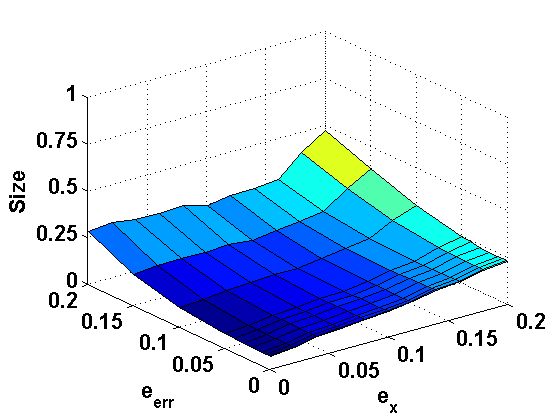}
		\label{FIG:9size_unknown_50_05}}
	~ 
	\subfloat[$n=50$, $\tau=\gamma=1$]{
		\includegraphics[width=0.3\textwidth]{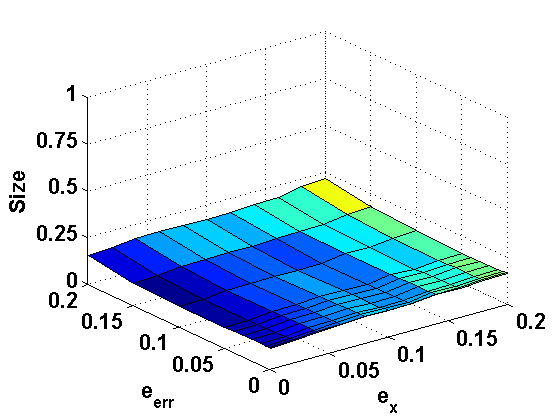}
		\label{FIG:9size_unknown_50_1}}
	\\
	\subfloat[$n=100$, $\tau=\gamma=0$]{
		\includegraphics[width=0.3\textwidth]{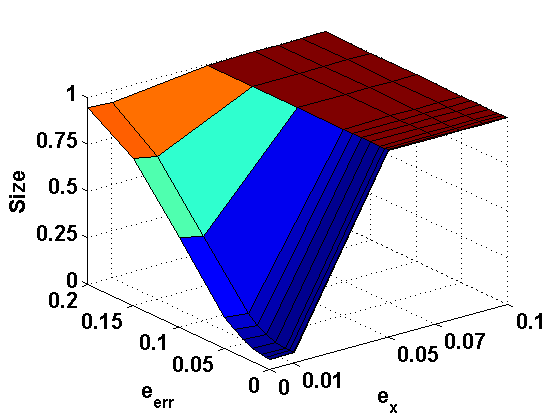}
		\label{FIG:9size_unknown_100_0}}
	~ 
	\subfloat[$n=100$, $\tau=\gamma=0.5$]{
		\includegraphics[width=0.3\textwidth]{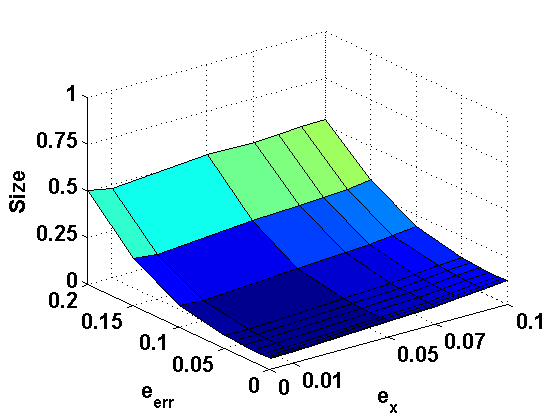}
		\label{FIG:9size_unknown_100_05}}
	~ 
	\subfloat[$n=100$, $\tau=\gamma=1$]{
		\includegraphics[width=0.3\textwidth]{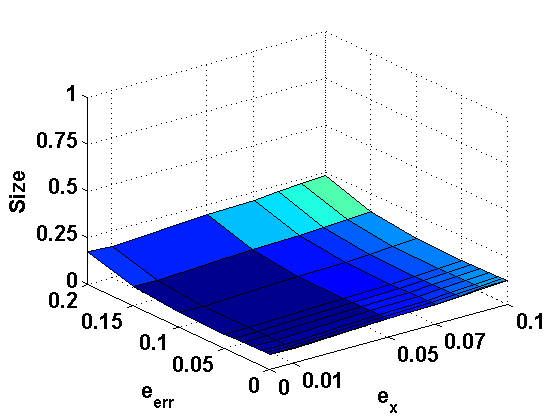}
		\label{FIG:9size_unknown_100_1}}
	\caption[Empirical size of the DPD based test of ${\boldsymbol{\beta}}$ with unknown $\sigma$]{Empirical size of the DPD based test of ${\boldsymbol{\beta}}$ with unknown $\sigma$ for different sample size $n$ and different $\tau=\gamma$}
	\label{FIG:9Size_unknown}
\end{figure}

\section{Empirical Illustrations}\label{SEC:9simulation_test}

\begin{figure}[h]
	\centering
	\subfloat[$n=30$, $\tau=\gamma=0$]{
		\includegraphics[width=0.3\textwidth]{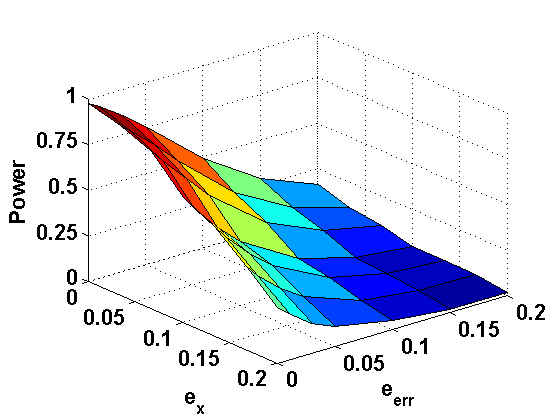}
		\label{FIG:9power_unknown_30_0}}
	~ 
	\subfloat[$n=30$, $\tau=\gamma=0.5$]{
		\includegraphics[width=0.3\textwidth]{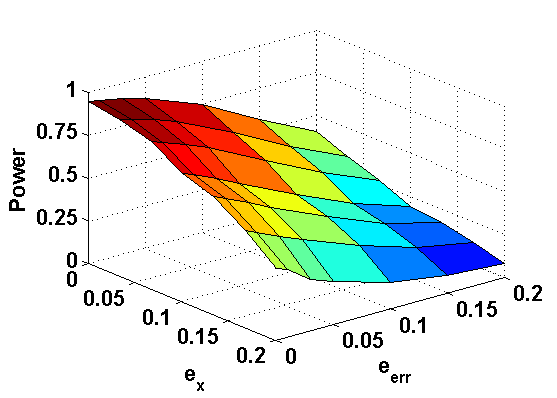}
		\label{FIG:9power_unknown_30_05}}
	~ 
	\subfloat[$n=30$, $\tau=\gamma=1$]{
		\includegraphics[width=0.3\textwidth]{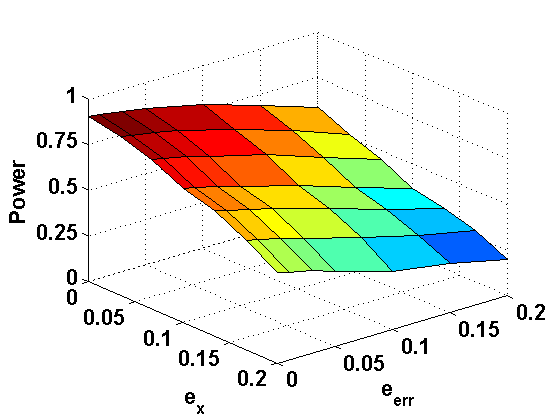}
		\label{FIG:9power_unknown_30_1}}
	\\
	\subfloat[$n=50$, $\tau=\gamma=0$]{
		\includegraphics[width=0.3\textwidth]{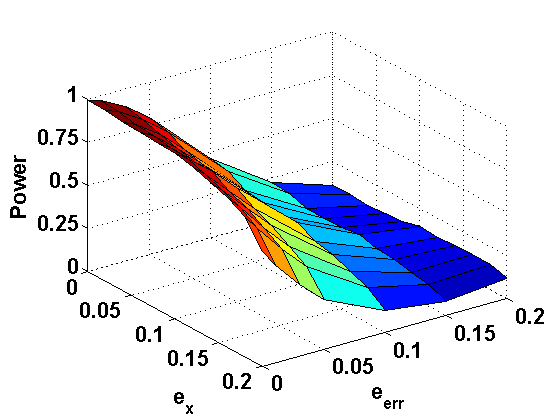}
		\label{FIG:9power_unknown_50_0}}
	~ 
	\subfloat[$n=50$, $\tau=\gamma=0.5$]{
		\includegraphics[width=0.3\textwidth]{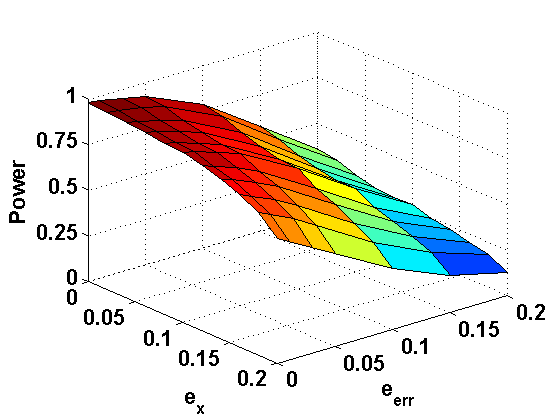}
		\label{FIG:9power_unknown_50_05}}
	~ 
	\subfloat[$n=50$, $\tau=\gamma=1$]{
		\includegraphics[width=0.3\textwidth]{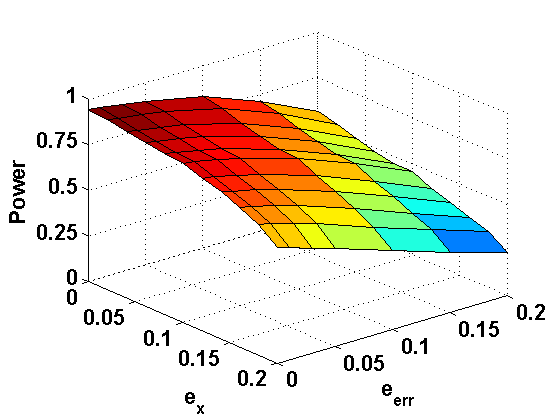}
		\label{FIG:9power_unknown_50_1}}
	\\
	\subfloat[$n=100$, $\tau=\gamma=0$]{
		\includegraphics[width=0.3\textwidth]{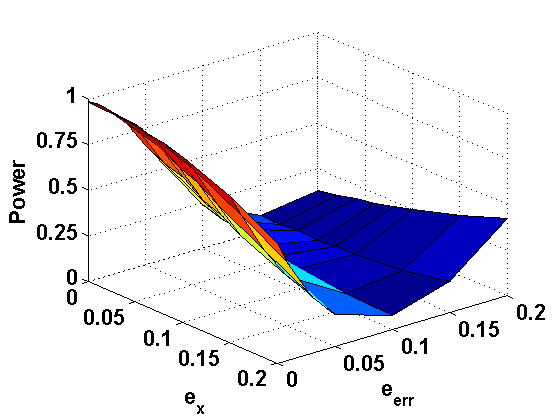}
		\label{FIG:9power_unknown_100_0}}
	~ 
	\subfloat[$n=100$, $\tau=\gamma=0.5$]{
		\includegraphics[width=0.3\textwidth]{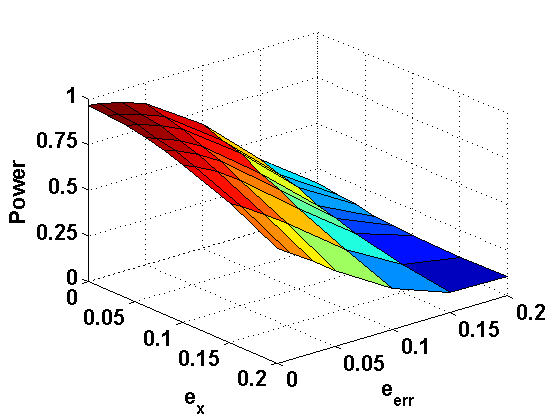}
		\label{FIG:9power_unknown_100_05}}
	~ 
	\subfloat[$n=100$, $\tau=\gamma=1$]{
		\includegraphics[width=0.3\textwidth]{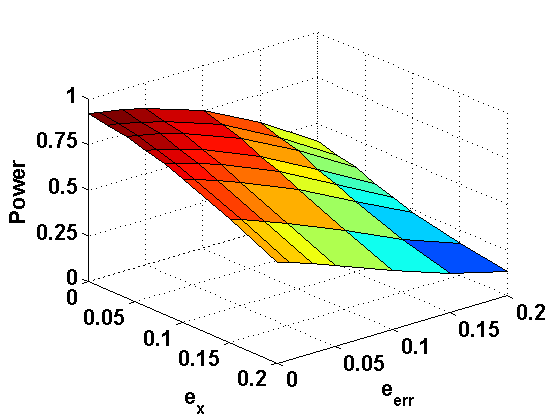}
		\label{FIG:9power_unknown_100_1}}
	\caption[Empirical power of the DPD based test of ${\boldsymbol{\beta}}$ with unknown $\sigma$]{Empirical power of the DPD based test of ${\boldsymbol{\beta}}$ with unknown $\sigma$ for different sample power $n$ and different  $\tau=\gamma$}
	\label{FIG:9power_unknown}
\end{figure}

We now illustrate the claimed robustness of the DPD based tests under an LRM with $\boldsymbol{x}_i = (1, z_i)^T$,
$z_i$ being fixed observation from $N(10, 5)$ distribution, and $\boldsymbol{\beta}=(\beta_1, \beta_2)^T$
for testing the composite null hypothesis $H_0: \boldsymbol{\beta}=(3, 2)^T$ assuming $\sigma$ unknown.
We replicate the simulation study of \cite{Ghosh/Basu:2017} which studied the robustness of simple null assuming $\sigma$ known.
We compute the empirical sizes and powers at the contiguous alternative $H_{1n}: \boldsymbol{\beta} = (3, 2)^T + \Delta_n$,
being $\Delta_n = \frac{1}{\sqrt{2n}}$, based on 1000 (independent)   LRM samples of sizes $n=30, 50$ and $100$.
In each sample, the errors are generated independently from $(1-e_{err}) N(0, 3) + e_{err} N(10, 3)$ distribution 
yielding $100e_{err}\%$ outliers in responses with true $\sigma=\sqrt{3}$.
We also simultaneously study the effect of leverage points;
randomly $100e_{x}\%$  of $z_i$s are replaced by observations from $N(16, 5)$ distribution
or by $\left[{{x}_i}(\frac{2-{{\Delta}}}{2})^2 - {{\Delta}}_n\right]$ respectively for size and power calculations.
These empirical sizes and powers are presented in Figures \ref{FIG:9Size_unknown} and 
\ref{FIG:9power_unknown} respectively for $\tau=\gamma=0$ (equivalent to LRT), $0.5$ and $1$.

Clearly the LRT ($\tau=\gamma=0$) is highly unstable with respect to both its size and power 
even for a fairly small contamination in either response or in design space. 
However, the DPD based tests with larger values of $\tau=\gamma$ are extremely robust 
against any kind of contamination in the data; their stability in both size and power increases as $\alpha$ increases.
This further justifies all theoretical robustness results derived here.

\vspace{-.4cm}
\section{Discussions}\label{SEC:conclusion}

This paper fills up the gaps of power and level robustness in the literature of 
DPD based robust tests for composite hypotheses.
The PIF and LIF are derived for general INH set-up and applied to the fixed-carrier linear regression model.
Further extension of the concept of PIF and LIF
for two or multi-sample problems under the INH set-up will be an interesting future work.

%

\appendix\small
\section{Notations}\label{APP:A}
\noindent
$\boldsymbol{0}_p$ denotes $p$-vector of zeros;
$\boldsymbol{O}_{p\times q} $ denotes null matrix of order $p\times q$
and $\boldsymbol{I}_{p} $ denotes identity matrix of order $p$.
\begin{eqnarray}
{\boldsymbol{u}_{i}}(y;{\boldsymbol{\theta}}) &=& \nabla \ln f_i(y;{\boldsymbol{\theta}})
~~\mbox{ with $\nabla$ representing gradient with respect to $\boldsymbol{\theta}$}\nonumber\\
H_n({\boldsymbol{\theta}}) &=& \frac{1}{n} \sum_{i=1}^n V_i(Y_i;{\boldsymbol{\theta}}), 
\mbox{ with }V_i(y;{\boldsymbol{\theta}}) =\int f_i(y;{\boldsymbol{\theta}})^{1+\tau} dy \nonumber
- \left(1+\frac{1}{\tau}\right) f_i(Y_i;{\boldsymbol{\theta}})^\tau\\
\boldsymbol\xi_i({\boldsymbol{\theta}}^g)  &=& 
\int {\boldsymbol{u}_{i}}(y;{\boldsymbol{\theta}}^g) f_{i}(y;{\boldsymbol{\theta}}^g)^\tau g_i(y)dy.
\label{EQ:xi}
\nonumber\\
\boldsymbol\Psi_n^\tau({\boldsymbol{\theta}}^g) &=& \frac{1}{n} \sum_{i=1}^n 
\left[\int {\boldsymbol{u}_{i}}(y;{\boldsymbol{\theta}}^g) {\boldsymbol{u}_{i}}^T(y;{\boldsymbol{\theta}}^g) 
f_i^{1+\tau}(y;{\boldsymbol{\theta}}^g) dy \right.\nonumber \\
&-& \left. \int \{ \mathbf{\nabla} {\boldsymbol{u}_{i}}(y;{\boldsymbol{\theta}}^g)  
+ \tau  {\boldsymbol{u}_{i}}(y;{\boldsymbol{\theta}}^g) {\boldsymbol{u}_{i}}^T(y;{\boldsymbol{\theta}}^g)\}
\{g_i(y) - f_{i}(y;{\boldsymbol{\theta}}^g)\} f_{i}(y;{\boldsymbol{\theta}}^g)^\tau  dy\right], \label{EQ:psi_n}
\nonumber\\
\boldsymbol\Omega_n^\tau({\boldsymbol{\theta}}^g) 
	&=& \frac{1}{n} \sum_{i=1}^n \left[ \int {\boldsymbol{u}_{i}}(y;{\boldsymbol{\theta}}^g) {\boldsymbol{u}_{i}}^T(y;{\boldsymbol{\theta}}^g) 
	f_{i}(y;{\boldsymbol{\theta}}^g)^{2\tau} g_i(y)dy 
	- \boldsymbol\xi_i({\boldsymbol{\theta}}^g)   \boldsymbol\xi_i^{T}({\boldsymbol{\theta}}^g)  \right],
	\label{EQ:omega_n}
\nonumber\\
\boldsymbol{P}_n^\tau({\boldsymbol{\theta}}) 
&=& \left[\frac{\nabla^2 H_n({\boldsymbol{\theta}})}{(1+\tau)}\right]^{-1} 
\left[\boldsymbol{I}_p - \boldsymbol\Upsilon({\boldsymbol{\theta}}) 
\boldsymbol\Upsilon^\ast({\boldsymbol{\theta}})^{-1} 
\boldsymbol\Upsilon({\boldsymbol{\theta}})^T [\nabla^2 H_n({\boldsymbol{\theta}})]^{-1} \right].
\nonumber\label{EQ:8P_n}\nonumber\\
\mbox{with } && \boldsymbol\Upsilon({\boldsymbol{\theta}}) 
= \frac{\partial \boldsymbol\upsilon({\boldsymbol{\theta}})}{\partial {\boldsymbol{\theta}}}
\mbox{ and }
\boldsymbol\Upsilon^\ast({\boldsymbol{\theta}}) = \left[\boldsymbol\Upsilon({\boldsymbol{\theta}})^T
[\nabla^2 H_n({\boldsymbol{\theta}})]^{-1}\boldsymbol\Upsilon({\boldsymbol{\theta}})\right]\nonumber\\
\boldsymbol{A}_n^\gamma({\boldsymbol{\theta}}) &=& 
\frac{1}{n} \sum_{i=1}^n \boldsymbol{A}_\gamma^{(i)}({\boldsymbol{\theta}})~~
\mbox{ with  }~~
\boldsymbol{A}_\gamma^{(i)}({\boldsymbol{\theta}}_0) = \nabla^2 d_\gamma(f_i(.;{\boldsymbol{\theta}}),f_i(.;{\boldsymbol{\theta}}_0))\big|_{{\boldsymbol{\theta}} = {\boldsymbol{\theta}}_0}.
\nonumber\\
\boldsymbol{M}_{j,\gamma}^{(i)}({\boldsymbol{\theta}}_1,{\boldsymbol{\theta}}_2) &=& 
\frac{\partial}{\partial\boldsymbol{\theta}_j}d_\gamma(f_i(.;{\boldsymbol{\theta}}_1),f_i(.;{\boldsymbol{\theta}}_2)),
~~j=1, 2; ~i=1, \ldots, n.\nonumber\\
\boldsymbol{A}_{j,k,\gamma}^{(i)}({\boldsymbol{\theta}}_1,{\boldsymbol{\theta}}_2) &=& 
\frac{\partial^2}{\partial\boldsymbol{\theta}_j\boldsymbol{\theta}_k}
d_\gamma(f_i(.;{\boldsymbol{\theta}}_1),f_i(.;{\boldsymbol{\theta}}_2)),
~~j, k=1, 2; ~i=1, \ldots, n.\nonumber
\end{eqnarray}

\end{document}